\documentclass[12pt]{article}

\usepackage{ifpdf}
\ifpdf
    \usepackage{graphicx}
 \else
    \usepackage[dvips]{graphicx}

\fi

\usepackage{amsmath,amscd,amssymb,amsthm,epsfig,array,hhline,longtable,float,eufrak}
\usepackage[alphabetic]{amsrefs}

\newcounter{assumption}[subsection]
\newenvironment{asm}%
  { %
    \begin{list}{\arabic{section}.\arabic{subsection}.\Alph{assumption}.}%
       { %
          \usecounter{assumption}%
\setlength{\leftmargin}{0pt}\setlength{\rightmargin}{0pt}%
            \setlength{\itemindent}{1.55cm}\setlength{\labelsep}{.2cm}%
            \setlength{\labelwidth}{1cm}%
       }%
   }%
  {\end{list}}

\newcounter{assumptionn}[subsection]
  { %
    \begin{list}{\arabic{section}.\arabic{subsection}.\arabic{assumptionn}.}%
       { %
          \usecounter{assumptionn}%
	\setlength{\topsep}{4pt}\setlength{\partopsep}{0pt}%
            \setlength{\itemsep}{2pt}%
            \setlength{\leftmargin}{0pt}\setlength{\rightmargin}{0pt}%
            \setlength{\itemindent}{1.55cm}\setlength{\labelsep}{.2cm}%
            \setlength{\labelwidth}{1cm}%
       }%
   }%
  {\end{list}}

\newcounter{assumptionnn}[subsubsection]
  { %
    \begin{list}{\arabic{section}.\arabic{subsection}.\arabic{subsubsection}.\Alph{assumptionnn}.}%
       { %
          \usecounter{assumptionnn}%
    \setlength{\topsep}{4pt}\setlength{\partopsep}{0pt}%
            \setlength{\itemsep}{2pt}%
            \setlength{\leftmargin}{0pt}\setlength{\rightmargin}{0pt}%
            \setlength{\itemindent}{1.55cm}\setlength{\labelsep}{.2cm}%
            \setlength{\labelwidth}{1cm}%
       }%
   }%
  {\end{list}}

\newcounter{assumptiona}[section]
  { %
    \begin{list}{\arabic{section}.\Alph{assumptiona}.}%
       { %
          \usecounter{assumptiona}%
\setlength{\leftmargin}{0pt}\setlength{\rightmargin}{0pt}%
            \setlength{\itemindent}{1.55cm}\setlength{\labelsep}{.2cm}%
            \setlength{\labelwidth}{1cm}%
       }%
   }%
  {\end{list}}

\theoremstyle{plain}

\newtheorem{thm}{Theorem}[section]

\newtheorem{lem}{Lemma}[section]
\newtheorem{prop}{Proposition}[section]
\newtheorem{corol}{Corollary}[section]

\theoremstyle{definition}

\newtheorem{defn}{Definition}[section]

\theoremstyle{remark}

\newtheorem{rmrk}{Remark}[section]

\numberwithin{equation}{section}

\renewcommand{\aa}{\alpha}
\newcommand{\bb}{\beta}
\renewcommand{\gg}{\gamma}
\newcommand{\GG}{\Gamma}
\newcommand{\dd}{\delta}
\newcommand{\DD}{\Delta}
\newcommand{\ee}{\varepsilon}

\renewcommand{\ll}{\lambda}
\newcommand{\LL}{\Lambda}
\renewcommand{\ss}{\sigma}
\newcommand{\ff}{\varphi}

\newcommand{\oo}{\omega}

\renewcommand{\SS}{\Sigma}

\renewcommand{\L}{\mathcal{L}}

\newcommand{\BB}{\mathbb{B}}
\newcommand{\C}{\mathbb{C}}

\newcommand{\PP}{\mathbb{P}}

\newcommand{\R}{\mathbb{R}}

\newcommand{\SP}{\mathbb{S}}

\newcommand{\AF}{\EuFrak{A}}

\newcommand{\SF}{\EuFrak{S}}

\newcommand{\<}{\langle}
\renewcommand{\>}{\rangle}

\newcommand{\qt}{\quad\text}

\newtheorem{example}{Example}[section]

\DeclareMathOperator{\Aut}{Aut}

\DeclareMathOperator{\GL}{GL}

\DeclareMathOperator{\Id}{Id}
\DeclareMathOperator{\AFF}{AFF}
\DeclareMathOperator{\Cor}{Cor}

\newcommand{\AAA}{\mathbb{A}}

\newcommand{\td}{\widetilde}
\renewcommand{\hat}{\widehat}
\renewcommand{\bar}{\overline}
\renewcommand{\tilde}{\widetilde}
\newcommand{\ga}{\mathfrak{g}}
\newcommand{\ha}{\mathfrak{h}}

\newcommand{\FB}{\mathbb{F}}

\begin{document}
\title{Enhanced Dynkin diagrams and Weyl orbits}
\author{E.\,B.\,Dynkin and A.\,N.\,Minchenko}
\maketitle

\begin{abstract}
The root system $\SS$ of a complex semisimple Lie algebra is uniquely determined by its basis (also called a simple root system). It is natural to ask whether all homomorphisms of root systems come from homomorphisms of their bases. Since  the Dynkin diagram of $\SS$ is, in general, not large enough to contain the diagrams of all subsystems of $\SS$, the answer to this question is negative. In this paper we introduce a canonical enlargement of a basis (called an enhanced basis) for which the stated question has a positive answer.  We use the name an enhanced Dynkin diagram for a diagram representing an enhanced basis.  These diagrams in combination with other new tools (mosets, core groups) allow to obtain a transparent picture of the natural partial order between Weyl orbits of subsystems in $\SS$. In this paper \textit{we consider only ADE root systems} (that is, systems represented by  simply laced Dynkin diagrams). The general case will be the subject of the next publication.
\end{abstract}

\tableofcontents


\section{Introduction}
A root system $\SS$ of a complex semisimple Lie algebra $\ga$ is a finite subset of a linear space $\ha(\R)^*$ dual to the real form $\ha(\R)$ of a Cartan subalgebra $\ha\subset\ga$. The Killing form allows to identify $\ha(\R)^*$ with $\ha(\R)$ and to interpret it as a Euclidean space. We say that $\td\SS\subset\SS$ is a subsystem of $\SS$ if it is a root system in its span and, for every $\aa,\bb\in\td\SS$, $\aa+\bb\in\SS$ implies $\aa+\bb\in\td\SS$. To every such $\td\SS$ there corresponds a semisimple subalgebra $\td\ga\subset\ga$ generated by root vectors $e_\aa\in\td\ga$ where $\aa\in\td\SS$. Subalgebras $\td\ga$, called \emph{regular}, play a fundamental role in the classification of all semisimple subalgebras of semisimple Lie algebras \cite{Dyn52}. There is a natural partial order between conjugacy classes of regular subalgebras of $\ga$. Namely, $C_1\prec C_2$ if $\ga_1\subset\ga_2$ for some $\ga_1\in C_1$ and $\ga_2\in C_2$. Enhanced Dynkin diagrams appeared as a result of our attempts to understand better this order.

A classification of complex semisimple Lie algebras is provided by the list of Dynkin diagrams. The Dynkin diagram $\GG$ of $\ga$ has nodes representing roots in a basis $\Pi\subset\SS$ and bonds describing relations between the roots. Subdiagrams of $\GG$ are Dynkin diagrams of  regular subalgebras of $\ga$. However, not all regular subalgebras can be obtained this way. Besides, non-conjugate regular subalgebras can have identical Dynkin diagrams. Both problems can be  efficiently solved by using  enhanced Dynkin diagrams.

In this paper we associate a diagram  with every subset $\LL$ of $\SS$. If $\LL$ is a basis, it coincides with the Dynkin diagram $\GG$ of $\ga$. We  construct an  enhancement  $\DD$ of $\GG$ by a recursive procedure which we call the completion. At each step, an extra node is introduced and connected by bonds with a proper part of already  introduced nodes. Enhanced Dynkin diagrams for simple $ADE$ algebras $\ga$ are presented on Figure \ref{Fig}.\footnote{Notation of roots is coordinated with that on Dynkin and extended Dynkin diagrams shown on Figure \ref{A}.}

Conjugacy classes of regular subalgebras of $\ga$ (with respect to the group of inner automorphisms of $\ga$)
correspond bijectively to conjugacy classes of subsystems of $\SS$ (with respect to the Weyl group $W$ of $\SS$). The correspondence preserves the natural partial order. Group $W$ acts on the set $\PP$ of bases of all subsystems $\td\SS\subset\SS$. Following \cite{Dyn52},  we use a name  \emph{$\Pi$-systems} for elements of $\PP$ . Every $\Pi$-system is contained in a unique subsystem of $\SS$. This implies a bijection between $W$-orbits in $\PP$ and  classes of conjugate subsystems of $\SS$.

A crucial step is a transaction from elements $\aa$ of $\SS$ to pairs $(\aa,-\aa)$. We call them \emph{projective roots}. To every root $\aa$ there corresponds a projective root $(\aa,-\aa)$, and to every set $\LL\subset\SS$ there corresponds a set of projective roots $p(\LL)$. We put $p(\SS)=S$ and we call $p(\LL)\subset S$ a projective basis if $\LL$ is a basis of $\SS$.  Projective subsystems and projective $\Pi$-systems of $S$ are defined in a similar way. Map $p$ induces a bijection between the class of all subsystems of $\SS$ and the class of all projective subsystems. It also induces a bijection between $W$-orbits in $p(\PP)$ and classes of conjugate subsystems of $\SS$.

If $L\subset S$ and if $|L|=k$,
\footnote{We denote by $|L|$ the cardinality of a set $L$.}
then $L=p(\LL)$ for  $2^k$ sets $\LL$ obtained from $L$ by  selecting one element from every pair $(\aa,-\aa)$. All $\LL$ have the same diagram which we assign as the diagram of $L$. We call a set of p rojective roots an \emph{enhanced basis}  if its diagram is isomorphic to the enhanced Dynkin diagram of  $\SS$. All enhanced bases  are conjugate. Besides, every projective $\Pi$-system is conjugate to a subset of an enhanced basis. We reduce the classification of all subsystems to the classification of orthogonal subsystems, i.e.\ those with bondless Dynkin diagrams. One of our tools is a Reduction Lemma which specifies in any projective $\Pi$-system $L_1$ elements $a$ with the property: if $f$ is an isometry of $L_1$ onto $L_2\subset S$ and if $f$ coincides with a $\tilde{w}\in W$, on $L_1\setminus\{a\}$, then $f=w$ on $L_1$ for some $w\in W$.
\footnote{In general, $w\neq\tilde w$.}

A special role belongs to maximal orthogonal subsets $M$ of $\SS$, which we call them \emph{mosets}. All mosets in $\SS$ are conjugate. To classify orthogonal subsystems in a  root system $\SS$, we consider a subgroup of $W$  formed by $w\in W$ stabilizing $M$ (we call it a \emph{core group}). We describe it in terms of the enhanced Dynkin diagram of $\SS$ and we apply  this result to characterize  conjugacy classes of  subsystems $\tilde\SS\subset\SS$ (and even homomorphisms $\tilde\SS\to\SS$ ) as well as a partial order between these classes. A byproduct of this investigation is an alternative way to obtain a known classification of regular subalgebras.

\subsection{Organization of the paper}
We deal with several classes of objects: $\AAA$ -- subsets of  a root system $\SS$; $\BB$ -- subsets of  a projective root system $S$;  $\C$ -- diagrams of $\AAA$ and $\BB$. The Weyl group $W$ acts on  $\AAA$, $\BB$ and $\C$. Our  goal is to investigate W-orbits in $\AAA$. Classes $\BB$ and $\C$ are tools to this end.

Root subsystems, bases (simple root systems), extended bases, $\Pi$-systems are subclasses of $\AAA$ commonly used in the theory of Lie algebras. All of them have shadows in $\BB$. If this can cause no confusion, we apply the same name and the same notation for each subclass and for its shadow. For instance, notation $\PP$ is used for the class of $\Pi$-systems in $\AAA$ and in $\BB$,
\footnote{However, to avoid confusion, we write sometimes $\PP(\SS)$ and $\PP(S)$.}
and a name $D_4$-sets is used for subsets of $\SS$ and for subsets of $S$ with  the Dynkin diagram of type $D_4$. We keep the same convention for  concepts like mosets and core groups introduced in this paper.

Often a statement formulated in the language of one of classes $\AAA, \BB, \C$
can be easily translated into the language of another class. We may state a proposition in a language of $\AAA$, prove it in terms of $\BB$ and apply it in the setting of $\C$.  This is done, for instance, in Section 4 where Classification Theorems \ref{thm_main(1)}, \ref{thm_main(2)}, \ref{thm_main(3)} are stated in terms of root systems, proved in terms of projective root systems and applied to description of Weyl orbits in terms of  enhanced Dynkin diagrams. In Section 2 completion is defined for subsets of $\SS$ but  constructed by using subdiagrams  of an enhanced Dynkin diagram.

In Section 1  we place the basic diagrams after Historical notes  because  references to them are spread over the  entire paper.

We start Section 2 with introducing diagrams $\GG$ of $\LL\subset \SS$ and $\DD$ of $L\subset S$. We define a class of complete subsets of $\SS$ and we introduce a completion $\bar X$ of $X\subset\SS$ as a minimal complete extension of $X$. If $X$ is a basis, then $\bar X$ is an enhanced basis and its diagram is an enhanced Dynkin diagram of $\SS$. We denote it $\DD(\SS)$. A
recursive  procedure for constructing a completion is illustrated by an example on Figure \ref{E7}.We conclude the section by proving important properties of enhanced and extended bases and automorphisms of their diagrams. In particular, we establish  that every $\Pi$-system is conjugate to a subset of an enhanced basis.

 A crucial role  in investigating which $\Pi$-subsystems of an enhanced basis are conjugate is played by Reduction Lemma that is proved at the beginning of Section 3. The rest of this  section is devoted to a description of mosets and core groups for all irreducible root systems.

In the Appendix we formulate basic facts on root systems, their bases and automorphisms and on the Weyl groups. For proofs we refer to the literature. At the end of the Appendix we list notation used in the paper.

\subsection{Historical notes}
An initial motivation for investigating subalgebras of Lie algebras was the problem of a classification of primitive transformation groups posed by S.~Lie and solved by him in spaces of 1, 2 and 3 dimensions. The case when the group is not semisimple was examined by V.~V.~Morozov in~1939.

Investigating semisimple subalgebras in arbitrary Lie algebras can be reduced to studying such subalgebras in semisimple algebras. For classical algebras $A_n$, $B_n$, $C_n$, $D_n$ semisimple subalgebras were studied by Malcev~\cite{Mal44}. The case of exceptional algebras $G_2$, $F_4$, $E_6$, $E_7$, $E_8$ was treated by Dynkin~\cite{Dyn52}. A principal step in this treatment was the classification of all regular semisimple subalgebras. For every simple Lie algebra, a list of elements $\LL_1,\dots,\LL_m$ of $\PP$ was given such that every Weyl orbit in $\PP$ is represented by  subsets of $\LL_k, k=1,\dots, m$. It was determined by direct case by case computations which of these subsets represent the same orbit. Only the results of these computations are presented in~\cite{Dyn52}.
The results of \cite{Dyn52} were used by a number of authors.

 Recently Oshima \cite{Osh07} verified them in his own way by introducing as a new tool the relation between orbits of $\LL$ and of its orthogonal complement in~$\SS$. His paper contains a number of useful tables.

Dynkin in \cite{Dyn52} applied his method also to subalgebras of classical algebras. However some time ago he discovered that his  description of orbits in $D_{2m}$ (more explicit than that in \cite{Mal44}) is incomplete. He also wanted to have a picture of the partial order between subalgebras which is very difficult to get from the classification in~\cite{Dyn52}. A discussion of these problems resulted in the beginning of the collaboration between  the authors of the present article. The ideas of introducing projective roots, enhanced Dynkin diagrams and mosets  should be credited to the second author.

Remarkably, slight modifications of enhanced Dynkin diagrams appeared recently in a paper of McKee and Smyth~\cite{McKS07}. Their diagrams served for a description of all integer symmetric matrices with all eigenvalues in the interval~$[-2, 2]$. Investigation of the relation between the two classes of diagrams may be a subject of the future work.

\subsection{Basic diagrams}
\begin{figure}[H]
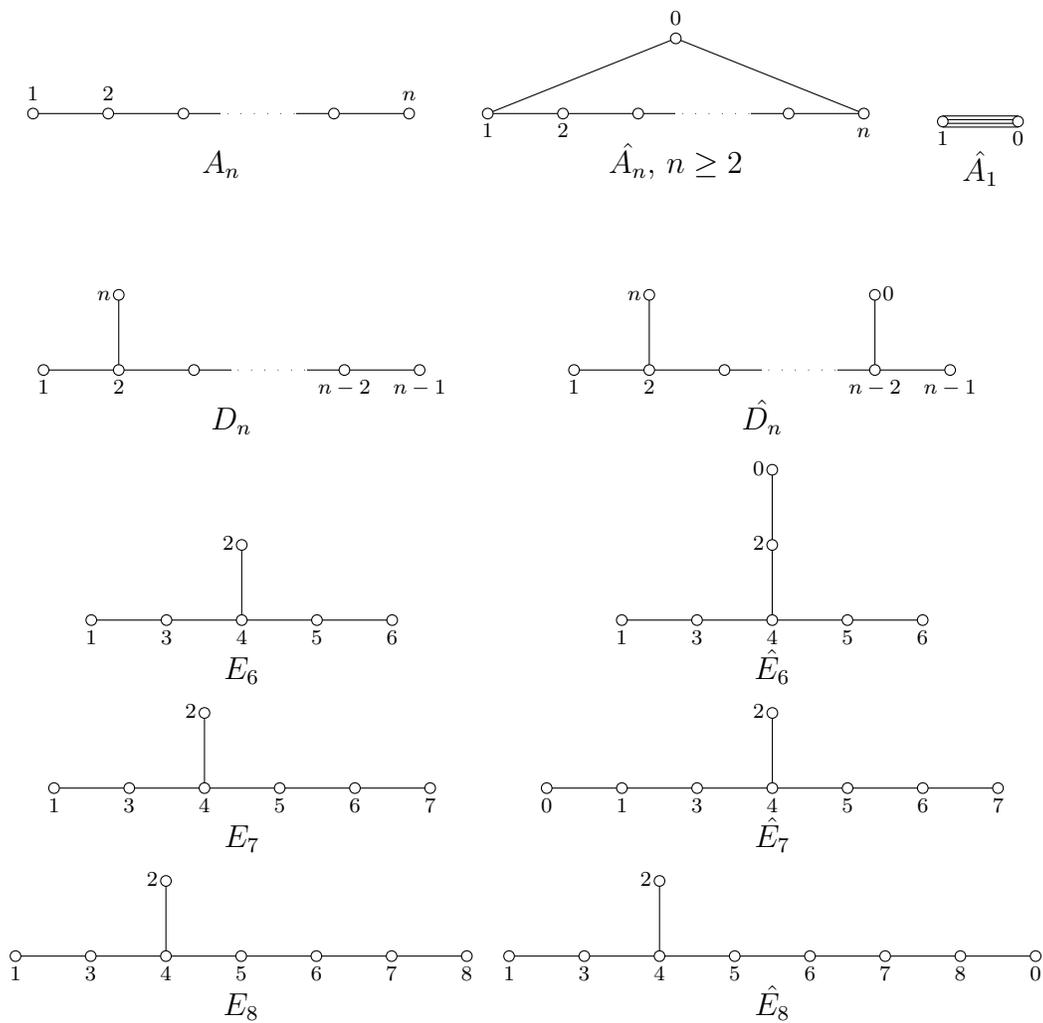

\begin{gather*}
\includegraphics{an.1}\qquad \includegraphics{extan.1}\qquad \includegraphics{exta1.1}\\
\end{gather*}
$$
\begin{array}{cc}
\includegraphics{dn.1}& \includegraphics{extdn.1}\\
\includegraphics{e6.1}& \includegraphics{exte6.1}\\
\includegraphics{e7.1}& \includegraphics{exte7.1}\\
\includegraphics{e8.1}& \includegraphics{exte8.1}
\end{array}
$$
\caption{Dynkin diagrams and extended Dynkin diagrams of $ADE$ systems}\label{A}
\end{figure}

\newpage
\begin{figure}[h]
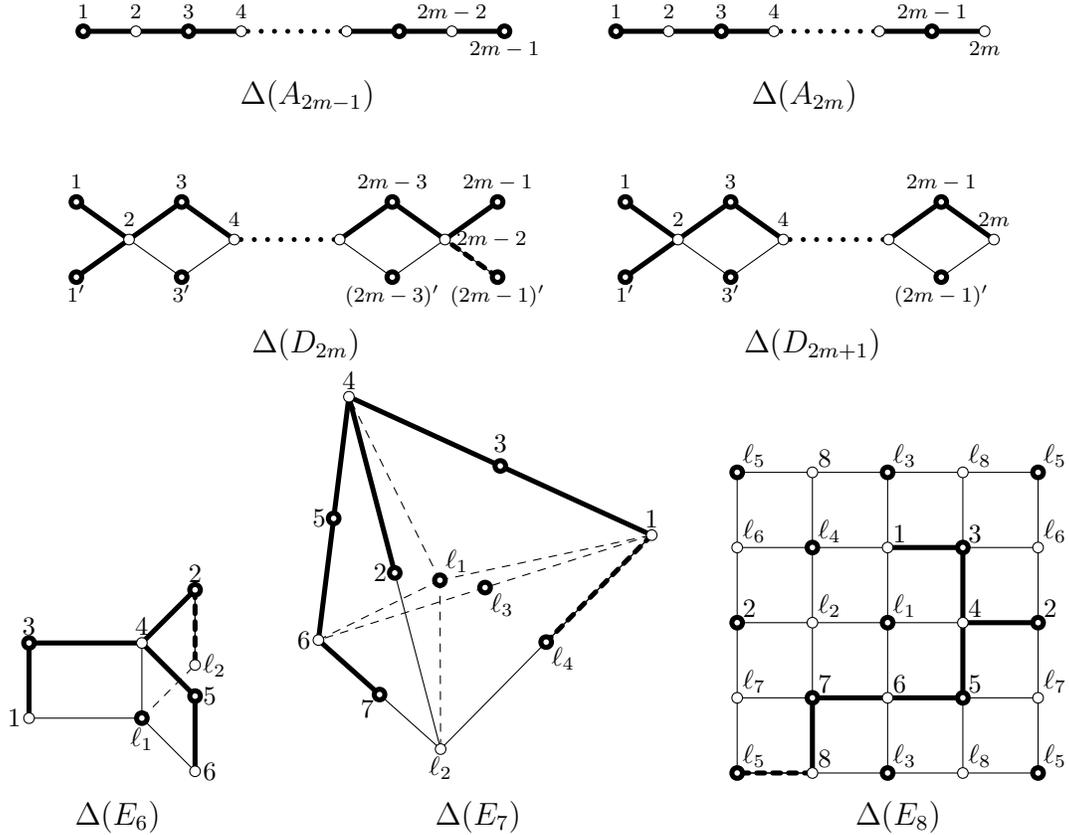

\begin{gather*}
\genfrac{}{}{0pt}{}{\includegraphics{anodd.1}}{\DD(A_{2m-1})}\qquad \genfrac{}{}{0pt}{}{\includegraphics{aneven.1}}{\DD(A_{2m})}\\
\\
\genfrac{}{}{0pt}{}{\includegraphics{dneven.1}}{\DD(D_{2m})}\qquad \genfrac{}{}{0pt}{}{\includegraphics{dnodd.1}}{\DD(D_{2m+1})}\\
 \genfrac{}{}{0pt}{}{\includegraphics{edde6.1}}{\DD(E_6)}\qquad \genfrac{}{}{0pt}{}{\includegraphics{edde7.1}}{\DD(E_7)}\qquad \genfrac{}{}{0pt}{}{\includegraphics{edde8.1}}{\DD(E_8)}\\
\end{gather*}
\caption{Enhanced diagrams for ADE root systems}\label{Fig}
\end{figure}

For every enhanced Dynkin diagram $\DD$ one of Dynkin diagram $\GG\subset \DD$ is marked by  boldfacing its bonds. Except $D_n$, we use for nodes of $\GG$  the same labels as on Figure~\ref{A}. The extra nodes for $E_6$, $E_7$ and $E_8$ are labeled $\ell_1,\ell_2,\dots$ in the order they are introduced by an algorithm described in Subsection \ref{subsection_completion}.  We draw separate diagrams for $D_{2m+1}$ and~$D_{2m}$. Label $1'$ is substituted for $n$ and extra nodes are labeled $3', 5',\dots,(2m-1)'$. For $E_6$ and $E_7$ it is convenient to embed the enhanced diagrams into the 3-dimensional space. In the case of $E_7$ the nodes are represented by the vertices of a tetrahedron, by the centers of its edges and by the center of the tetrahedron. Dashed lines are used for bonds shaded by faces. The enhanced diagram for $E_8$ is a $4\times 4$ lattice on the torus: each bottom node represents the same root as the corresponding top node and a similar identification is made for the left and right sides of the square.

The diagram $\DD(E_8)$ can be also interpreted as vertices and edges of a hypercube in $\R^4$.  $\DD(E_7)$ can be obtained by eliminating any vertex and all edges containing it. To get $\DD(E_6)$ one needs to drop any edge  and all edges on 2-faces that contain this edge.

Boldfaced nodes represent a maximal orthogonal subsystem of $S$.

\section{Enhanced Dynkin diagrams}
Starting from this point we consider only ADE root systems $\SS$. In this case the Cartan integers $(\aa|\bb)=\frac{2(\aa,\bb)}{(\bb,\bb)}$ are equal to $0$, $1$ or $-1$ for all non-collinear $\aa,\bb\in\SS$.

\subsection{Diagrams and projective diagrams}
We define a \emph{diagram $\GG_\LL$} of $\LL\subset\SS$  by assigning a node to every $\aa\in\LL$ and by connecting two nodes representing roots $\aa,\bb\in\LL$  by a single bond, if $(\aa,\bb)\neq 0$ and $\aa\neq -\bb$, and by a quadruple bond, if $\aa=-\bb$.\footnote{The number of bonds is equal to the integer $(\aa|\bb)(\bb|\aa)$ which can be only 0,1 or 4 in case of ADE root systems.} We also introduce   a \emph{projective diagram} $\DD_\LL$ of $\LL$ which is  obtained from $\GG_\LL$ by identifying two nodes corresponding to every pair $\aa,-\aa$ (and by dropping the quadruple bond connecting the nodes). Two nodes of $\DD_\LL$ are connected by a bond if and only if  they represent non-orthogonal roots. We call such nodes (and roots)  \emph{neighbors}. Subdiagrams of any diagram are in a 1-1 correspondence with subsets of nodes.

A subset $\LL\subset\SS$ is called \emph{symmetric}, if $\LL=-\LL$. There is a 1-1 correspondence between symmetric subsets of $\SS$ and all subsets of  the projective root system $S$.[The latter correspond bijectively to subdiagrams of the projective diagram of $\SS$.]

A set of roots is called \emph{irreducible} if it's diagram is connected. Every set $\LL\subset\SS$ is the union of its maximal irreducible subsets. We  call them \emph{irreducible components} of $\LL$.

A linearly independent subset $\LL$ of $\SS$ is a basis of a subsystem $\td\SS$ if and only if  $\aa-\bb\notin\SS$ for all $\aa,\bb\in\LL$. Recall that we use the name \emph{$\Pi$-systems} for  such $\LL$ and that we denote the class of all $\Pi$-systems by $\PP$.

A partial order in $\SS$ is associated with a basis $\Pi$ by the condition: $\aa\prec\bb$ if $\bb-\aa$ can be represented as the sum of elements of $\Pi$. If $\SS$ is irreducible, then there is a unique minimal root
$\aa\in\SS$ with respect to this order (\cite[Chapter VI, Proposition 25(i)]{Bou02}). The set
\[
\hat\Pi=\Pi\cup\{\aa\},
\]
is called an \emph{extended basis}. The diagram $\hat\GG=\GG_{\hat\Pi}$ is the \emph{extended Dynkin diagram} of $\SS$.
Similarly, we define \emph{extended $\Pi$-system} $\hat\LL$ for every irreducible $\Pi$-system $\LL$. 

We say that nodes $a_1,a_2,\dots,a_k$, $k\ge 3$  of a diagram form a cycle if  $a_i$ is a neighbor of $a_{i+1}$ for $1\le i\le k-1$, and $a_k$ is a neighbor of  $a_1$. A diagram  is \emph{acyclic}, if it contains no cycles.
 We say that a root $\ee$ is an \emph{end} of $\LL\subset\SS$ if it has no more than one neighbor. Every root set with an acyclic diagram contains an end.

\begin{prop}
\label{section}
Let $\Phi$ be a symmetric subset of  $\SS$ such that $\DD_\Phi$ is connected and acyclic. Then $\Phi=\LL\cup -\LL$ where $\LL$ is either a $\Pi$-system or an extended $\Pi$-system.
\end{prop}

\begin{proof}
Choose any $\LL\subset \Phi$ such that $\LL\cap -\LL=\emptyset$ and $\LL\cup -\LL=\Phi$. Since $\GG_\LL=\DD_\Phi$ is acyclic, we can assume (by replacing $\aa$ by $-\aa$ for some $\aa\in\LL$) that  $(\aa,\bb)\le 0$ for all $\aa,\bb\in \LL$. Now Proposition \ref{section} follows immediately from \cite[Lemma 5.1]{Dyn52}.
\end{proof}

Proposition \ref{section} implies that all connected acyclic diagrams are presented on Figure~\ref{A}. 

\subsection{Completion}\label{subsection_completion}
\begin{defn}
A symmetric subset $\Phi\subset\SS$ is called \emph{complete} if  the condition $\LL\subset\Phi$ implies $\hat\LL\subset \Phi$ for every $D_4$-set  $\LL\in\PP$ of type . The \emph{completion} $\bar X$ of  $X\subset\SS$ is the intersection of all complete subsets of $\SS$ containing $X$.
\end{defn}
Clearly,  $\bar X\subset\bar Y$ if  $X\subset Y$.

To construct $\bar X$ we introduce a sequence of symmetric sets
\[
\Phi_0\subset\Phi_1\subset\dots\Phi_i\dots\subset\Phi_k
\]
where $\Phi_0=X$, $\Phi_k=\bar X$ and, for every  incomplete $\Phi_{i-1}$, a diagram $\DD_{\Phi_i}$ is defined as an elementary extension of  $\DD_{\Phi_{i-1}}$.

An elementary extension $\DD'$ of any  diagram  $\DD$ is obtained by extending a $D_4$-subdiagram $D\subset\DD$ to a $\hat D_4$-subdiagram of $\DD_{\SS}$.  $\DD'$  contains an extra node $\ell$ connected by a bond with the branching node of $D$ and with every node of $\DD\setminus{D}$ that has 1 or 3 neighbors among the ends of $D$. The procedure is motivated by the following proposition.

\begin{prop}
\label{t1}
Let $\LL$ be a $\hat D_4$-subset of a symmetric set $\Phi$. If $\bb\in\SS\setminus\LL$ and if $\LL\cup\{\bb\}$ is irreducible, then $\bb$ has two  or four neighbors among end elements of~$\LL$.
\end{prop}

\begin{proof}
If $\LL=\{\aa_0,\aa_1,\aa_2,\aa_3,\gg\}$ with  $\aa_0,\aa_1,\aa_2,\aa_3$ representing the end nodes, then
$\aa_0+\aa_1+\aa_2+\aa_3+2\gg=0.$
\footnote{To prove this, evaluate the squared length of the  left side.}
This implies
\[
(\aa_0|\bb)+(\aa_1|\bb)+(\aa_2|\bb)+(\aa_3|\bb)+2(\gg|\bb)=0.
\]
Taking this relation modulo 2, we obtain that $\bb$ is the neighbor of an even number of the ends in $\LL$.
\end{proof}

\begin{figure}[H]
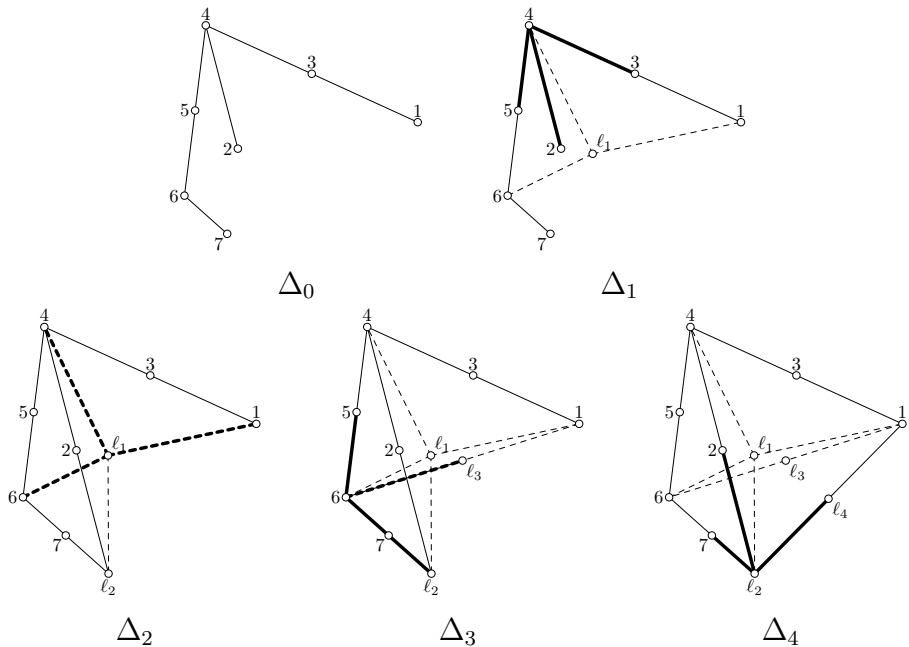

\begin{gather*}
 \genfrac{}{}{0pt}{}{\includegraphics[scale=.7]{e70.1}}{\DD_0}\qquad \genfrac{}{}{0pt}{}{\includegraphics[scale=.7]{e71.1}}{\DD_1}\\
 \genfrac{}{}{0pt}{}{\includegraphics[scale=.7]{e72.1}}{\DD_2}\qquad \genfrac{}{}{0pt}{}{\includegraphics[scale=.7]{e73.1}}{\DD_3}\qquad
 \genfrac{}{}{0pt}{}{\includegraphics[scale=.7]{e74.1}}{\DD_4}\\
\end{gather*}
\caption{Detail for $E_7$}\label{E7}
\end{figure}

On Figure \ref{E7} we illustrate this procedure on an example of a symmetric set $X$ described by the Dynkin diagram $\DD_0$ of $E_7$. We start with  its $D_4$-subdiagram $\{5,2,3,4\}$ (in notation of Figure \ref{A}).  We get a diagram $\DD_1$ by introducing an extra node $\ell_1$ connected with nodes 1, 4 and~6.  In $\DD_1$, we choose a $D_4$-subdiagram $\{1,4,6,\ell_1\}$
\footnote{Another possible choice is $\{\ell_1,5,7,6\}$.}
 and we introduce an extra node $\ell_2$ connected with $\ell_1$, 2 and~7. The result of the second step is $\DD_2$. We arrive at $\DD_3$ by extending a subdiagram $\{4,\ell_1,\ell_2,6\}$ of $\DD_2$ by a node $\ell_3$. Finally, we get a  diagram $\DD_4$ by introducing $\ell_4$ connected with $\ell_2$ and~1.  $\DD_4$ represents $\bar X$.
(A different selection of elementary extensions would lead to an isomorphic diagram with  a permutation of labels $\ell_1,\ell_2,\ell_3,\ell_4$.)

\subsection{Enhanced bases and enhanced Dynkin diagrams}
\begin{defn}
An \emph{enhanced basis} of $\SS$ is the completion of a basis of $\SS$. Its projective diagram is called an \emph{enhanced Dynkin diagram} of $\SS$.
\end{defn}

Conjugacy of all bases implies conjugacy of all enhanced bases, hence isomorphism of all enhanced Dynkin diagrams of $\SS$. These diagrams for irreducible ADE root systems are presented on Figure \ref{Fig}.

A fundamental  property of enhanced bases is established in the following theorem.
\begin{thm}\label{thm_rich}
Any $\Pi$-system is conjugate to a subset of an enhanced basis.
\end{thm}

In the proof of this property we use   a class of transformations of $\PP$ introduced in \cite{Dyn52}.
\begin{defn}
 An \textit{elementary transformation} of an irreducible $\Pi$-system $\LL$ is a set  $\LL^\aa=\hat\LL\setminus\{\aa\}$, where $\aa\in\hat\LL$. An elementary transformation of an arbitrary $\Pi$-system $\LL$  is an elementary transformation of one of its irreducible components. We call this transformation \emph{trivial} if  $\LL^\aa$ is isomorphic to $\LL$.\footnotemark
\end{defn}
\footnotetext{A bijection $\ff:\LL_1\to\LL_2$ is an isomorphism if it preserves the bracket $(\cdot|\cdot)$.}

It follows from \cite[Theorems 5.2, 5.3]{Dyn52} that:
\begin{prop}
\label{elem_trans}
Let $\SS$ be a root system of rank $n$. Every $\LL\in\PP$ with $|\LL|=n$   can be obtained from any basis of $\SS$ by a sequence of elementary transformations. If  $\SS$ is a  system of rank $n$, then  any $\LL\in\PP$ with $|\LL|< n$ is contained in a $\Pi$-system with $n$ elements.
\end{prop}

Proposition \ref{Phi} implies:
\begin{prop}
\label{trivial}
All trivial elementary transformations of $\LL\in\PP$ are conjugate.
\end{prop}

\begin{prop}
\label{dkdl_conj}
Let $\SS=D_n$\footnote{We abuse notation by writing, e.g.\ $X=D_n$ if $X$ is a root system, projective root system or a $\Pi$-system of type $D_n$.} and $2\le k\le n-2$.\footnote{We have $D_2=A_1+A_1$.} Then all $\Pi$-systems of type $D_k+D_{n-k}$ in $\SS$ are conjugate.
\end{prop}

\begin{proof}
Fix a basis $\LL$. Note that if $\LL^\aa$ has type $D_k+D_{n-k}$, then $\aa$ corresponds to either label $k$ or to $n-k$ on Figure \ref{A}. Clearly, these two $\Pi$-systems are conjugate by a $w\in W_{\hat{\Pi}}$ moving node $1$ to $n-1$. [Such $w$ exists by Proposition ~\ref{Phi}.] Finally, note that a $\Pi$-system of type $D_k+D_{n-k}$ cannot be obtained from $\LL$ by a sequence of more than one non-trivial elementary transformations, since the number of irreducible components of the transformation increases after each step.
\end{proof}

\begin{proof}[Proof of Theorem \ref{thm_rich}]
Let $\Phi\subset\SS$ be an enhanced basis.
We need to demonstrate  that $\PP$ is  contained in the class $\L$ of  $\LL$ with the property: $w(\LL)\subset\Phi$ for some $w\in W$. Class $\L$ is preserved under the action of $W$ and it contains, with every $\LL$, all its subsets. Since $\Phi$ is the completion of a basis, we need only to demonstrate that  $\L$ contains, with every $\LL\in\PP$, all elementary transformations of~$\LL$. It is sufficient to check this for irreducible $\LL$.
 For $\LL=D_{2m}$, $E_6$, $E_7$ or $E_8$, this is true because $\Phi$ contains $\hat{\LL}$ (see Figure~\ref{Fig}). For type $A_n$ all elementary transformations are trivial,  hence they belong to $\L$ by Proposition \ref{trivial}. Finally, for type  $D_{2m+1}$ all elementary transformations   are of type $D_{k_1}+D_{k_2}$ with $k_1+k_2=2m+1$, and it is clear from Proposition \ref{dkdl_conj} and Figure~\ref{Fig} that they belong to~$\L$.
\end{proof}

For any group $G$ of transformations of $\SS$ and any  $X\subset\SS$, we put $G_X=\{g\in G: g(X)=X\}$. In particular, $W_X$ stands for the subgroup of the Weyl group preserving $X$.

\begin{prop}
\label{conj_inside}
All bases of $\SS$ contained in an enhanced basis $\Phi$ are conjugate by $w\in W_\Phi$.
\end{prop}
\begin{proof}
Suppose $\Phi=\bar\Pi$, where $\Pi\subset\SS$ is a basis. For any other basis $\Pi'\subset\Phi$, by Proposition \ref{base_conj},  there is a $w\in W$ such that $w(\Pi)=\Pi'$ and therefore $w(\bar\Pi)=\bar{\Pi'}$.  Since a complete $\Phi$ contains $\Pi'$, it contains $\bar{\Pi'}=w(\Phi)$, and $w(\Phi)=\Phi$ because  $|w(\Phi)|=|\Phi|$.
\end{proof}


\begin{prop}
\label{aut_enhanced}
Let $\Pi$ be a basis of $\SS$ and $\Phi=\Pi$, $\hat\Pi$ or $\bar\Pi$. Every automorphism of $\DD_\Phi$ is induced by an automorphism  of $\SS$. The latter can be chosen in $W$, if  $\SS$ is  not of type $D_{2n}$.
\end{prop}

\begin{proof}
Let $f$ be an automorphism of  $\DD=\DD_\Phi$. The Dynkin diagrams $\GG=\GG_\Pi$ and $\GG'=f(\GG)$ are subdiagrams of  $\DD$, and $\GG'=\GG_{\Pi'}$ for  a basis  $\Pi'\subset\Phi$. By Proposition \ref{conj_inside}, $w(\Pi')=\Pi$ for some $w\in W_\Phi$ and therefore $w(\GG')=\GG$ and $wf(\GG)=\GG$. By Proposition \ref{automorphisms}, the map $wf$ is induced by an automorphism $\ff$ of $\SS$. Automorphism $\psi=\ff^{-1}wf$ of $\DD$ is trivial on $\GG$. This implies cases $\Phi=\Pi$ and $\Phi=\hat\Pi$. If  $\Phi=\bar\Pi$, then $\DD$ is the enhanced  Dynkin diagram of $\SS$. There exists a sequence of elementary extensions $\GG=\DD_0\subset\DD_1\subset\cdots\subset\DD_k=\DD$. We prove by induction that  $\psi$  is trivial on $\DD_i$ for $i=0,\dots, k$. Therefore  $f$ is induced by $w^{-1}\ff\in\Aut\SS$.

If $\SS=D_{2n}$, then, by Proposition \ref{involution}, $\Aut\SS$ is generated by $W$ and $-\Id$. Since $-\Id$ acts trivially on $\DD$, $\Aut\SS$ and $W$ induce the same group of transformations on $\DD$.
\end{proof}

\section{Mosets and core groups}
\subsection{Reduction Lemma}
We consider a projective root system $S$ and we put
\[\<a|b\>=\left|(\aa|\bb)\right|,
\]
for every $a=(\aa,-\aa), b=(\bb,-\bb)\in S$. By replacing $(\aa|\bb)$ by $\<a|b\>$ we introduce  the concepts of orthogonality, irreducibility and
 isomorphism for subsets of $S$.

We denote $I(\LL,\SS)$ the set of all embeddings of $\LL\subset\SS$ into $\SS$ preserving $(\cdot|\cdot)$.\footnote{We will call elements of $I(\LL,\SS$ embeddings.} For a $f\in I(\LL,\SS)$ we write $f\in W(\LL,\SS)$ if there exists a $w\in W$ such that $f=w$ on $\LL$. Clearly, $W(\LL,\SS)\subset I(\LL,\SS)$. In a similar way we introduce classes $W(L,S)\subset I(L,S)$.

\begin{lem}[Reduction Lemma]
\label{lem_red}
Suppose $f\in I(L,S)$ and $a$ is the neighbor of  an end $e$ of $L$. Then
$f\in W(L,S)$ if $\td f\in W(\td L,S)$, where $\td f$ is the restriction of $f$ to $\td L=L\setminus\{a\}$.
\footnote{This Lemma is false for $S$ substituted by $\SS$. For instance, suppose $\SS$ is of type $A_2$, $\LL$ is its basis and $f=-\Id$. Any $\aa\in\LL$ is the neighbor of an end and the restriction $\td f$ of $f$ to $\td \LL=\LL\setminus\{\aa\}$ belongs to $W(\tilde\LL,\SS)$ but $f\notin W(\LL,\SS)$.}
\end{lem}

Note that reflections
\[
s_\aa(\bb)=\bb-(\bb|\aa)\aa,\quad\bb\in\SS,
\]
satisfy the condition
\[
s_\aa s_\bb s_\aa=s_\gg\qt{where }\gg=s_\aa(\bb), \qt{for every }\aa,\bb\in\SS.
\]

\begin{lem}
\label{lem_red1}
Let $\td L\subset L\subset S$, $f\in I(L,S)$, $f(a)=a$ for all $a\in \td L$. If $a\in L$ and $\<a|f(a)\>\neq 0$, then the restriction $\td f$ of $f$ to $\td L\cup\{a\}$ belongs to $W(\td L\cup\{a\},S)$.
\end{lem}
\begin{proof}
Set $b=f(a)$. If $a=b$, then $\td f=\Id$. Suppose $a\neq b$. Then there exist $\aa,\bb\in\SS$ such that $p(\aa)=a$, $p(\bb)=b$ and $(\aa,\bb)<0$. Since $(\aa|\bb)=(\bb|\aa)=-1$, we have
$$
s_\aa(\bb)=s_\bb(\aa)=\aa+\bb=\gg\in\Sigma
$$
and $s_\gg=s_\aa s_\bb s_\aa$.

Let us show that $s_\gg(a)=b$ and $s_\gg(\ell)=\ell$ for all $\ell\in\td{L}$. We have
$$
s_\gg(\aa)=s_\aa s_\bb s_\aa(\aa)=-s_\aa s_\bb (\aa)=-s_\aa s_\aa (\bb)=-\bb.
$$
 Hence $s_\gg(a)=b$. Let $\ell\in\tilde{L}$. There exists $\ll\in\SS$ such that $p(\ll)=\ell$ and $(\aa,\ll)\le 0$. Since $\left|(\aa|\ll)\right|=\left|(\bb|\ll)\right|$
we get $(\gg|\ll)=(\aa|\ll)+(\bb|\ll)=0\,\text{or}\,-2.$ If $(\gg|\ll)=0$, then $s_\gg(\ll)=\ll$. If $(\gg|\ll)=-2$, then $\gg=-\ll$, $s_\gg(\ll)=-\ll$ and thus $s_\gg(\ell)=\ell$.
\end{proof}

\begin{proof}[Proof of Lemma \ref{lem_red}]
 We may assume that $\td{f}=\Id$. Suppose
$ \td f\in W( \td L,S)$. If $\<a|f(a)\>\neq 0$, then  $ f\in W( L,S)$ by Lemma \ref{lem_red1}.

Now suppose $\<a|f(a)\>=0$. Let $\ee\in\SS$ be such that $p(\ee)=e$. Set $f'=s_\ee f\in I(L,S)$ and let $\ell\in\td L\setminus\{e\}$.  We have  $f'(\ell)=\ell$ because $\<e|\ell\>=0$. Since $\<a|f'(a)\>\neq 0$, $f'\in W(\td L,S)$ by Lemma \ref{lem_red1}. Hence $f=s_\ee f'\in W(\td L,S)$.
\end{proof}

\subsection{Mosets}
\begin{defn}
An orthogonal set $M\subset X\subset\SS$ is a \emph{moset} in $X$ if none element of $X\setminus M$ is orthogonal to $M$.
 A moset $O$ of an irreducible $X$ is  \emph{perfect} if $X\setminus O$ is orthogonal and $|O|\ge |X\setminus O|$. A perfect moset of an arbitrary $X$ is the union of perfect mosets of irreducible components of $X$.
\end{defn}

Every $\Pi$-system has a perfect moset. Indeed,  every basis $\Pi$ of an irreducible root system is the union  of two disjoint  orthogonal subsets and at least one of them is  a perfect moset in $\Pi$.

Mosets and perfect mosets in $S$ are defined similarly. Clearly,
$p(M)$ is a moset in $S$ if  $M$ is a moset in $\SS$.

\begin{thm}
\label{thm_perfect}
Let $O$ be a perfect moset in $L\in\PP(S)$ and let $f\in I(L,S)$. Then
$f\in W(L,S)$ if $\td f=f|_O\in W(O,S)$.
\end{thm}
\begin{proof}
Dynkin diagrams on Figure \ref{A} demonstrate that any perfect moset in $L$ contains at least one of the end nodes $e$. Hence it does not contain the neighbor $a$ of $e$, and we can delete $a$. As the result we get another perfect moset of a $\Pi$-system. Therefore we get Theorem \ref{thm_perfect} by using an induction in $|L\setminus{O}|$ and by applying  Lemma \ref{lem_red}.
\end{proof}

\begin{asm}
\item\label{1.8.a}
$M\subset\SS$ is a moset if and only if its intersection with every irreducible component of $\SS$ is a moset in this component.

\smallskip
The set $\Psi_X(\SS)\subset\SS$ of all roots orthogonal to $X\subset\SS$ forms a subsystem in $\SS$. For every $\LL\supset X$, we denote by $\Psi_X(\LL)$ the set of all $\bb\in\LL$  orthogonal to $X$.

\item\label{1.8.A}
If $\aa\in M\subset\SS$, then $M$ is a moset of $\SS$ if and only if $M\setminus\{\aa\}$ is a moset of $\Psi_\aa(\SS)$.

\smallskip

\begin{prop}
\label{moset_conj}
All mosets in $\SS$ are conjugate and therefore they have the same cardinality $\mu(\SS)$.
\end{prop}

This follows by induction in rank of $\SS$ from Proposition \ref{root_conj} and~\ref{1.8.a},~\ref{1.8.A}.

\smallskip
Since every orthogonal set is contained in a moset, we have
\item\label{1.8.C}
An orthogonal set with $m$  elements  is a moset if and only if $m=\mu(\SS)$.

\smallskip
It follows from \ref{1.8.A} that
\item\label{1.8.D}
$\mu(\SS)=\mu(\Psi_\aa(\SS))+1$.

\smallskip
\item\label{1.8.E}
Values of $\mu(\SS)$ for irreducible $\SS$ are given by the table:

\begin{longtable}{|c|c|c|c|c|c|}
\caption{Cardinality of mosets} \label{mu}\\
\hline
$\SS:$ & $A_n$ & $D_n$ & $E_6$ & $E_7$ & $E_8$ \\
\hhline{|=|=|=|=|=|=|}
$\mu(\SS):$ & $q_{n+1}$ & $2q_n$ & $4$ & $7$ & $8$ \\
\hline
\end{longtable}
where $q_n$ is a maximal integer not exceeding $\frac{n}{2}$.

\smallskip
This is an implication of \ref{1.8.D} and the following result.

\item\label{1.4.I}
Types of $\Psi_\aa(\SS)$ for irreducible $\SS$ are given in Table \ref{tocl}.

\begin{longtable}{|c|c|c|c|c|c|}
\caption{Types of the orthogonal complements to roots}\label{tocl}\\
\hline
$\SS:$ & $A_n$ & $D_n$ & $E_6$ & $E_7$ & $E_8$\\
\hhline{|=|=|=|=|=|=|}
$\Psi_\aa(\SS):$ & $A_{n-2}$ & $D_{n-2}+A_1$ & $A_5$ & $D_6$ & $E_7$\\
\hline
\end{longtable}

[Since all roots are conjugate, their orthogonal complements have the same type.]

\begin{proof}
 Let  $\aa$ be the minimal root with respect to a basis $\Pi$ of $\SS$. For every $0<\bb\in\SS$ we have $$\bb=\sum_{\gg\in\Pi} k_\gg \gg,$$ where $k_\gg\ge 0$ for all $\gg\in\Pi$. Therefore $\bb\in\Psi_\aa(\SS)$ if and only if $\sum_{\gg\in\Pi} k_\gg(\aa,\gg)=0$. Since $\aa$ is minimal, $(\aa,\gg)\le 0$ for all $\gg\in\Pi$. Hence  $\bb\in\Psi_\aa(\SS)$ if and only if $k_\gg=0$ for all $\gg$ such that $(\aa,\gg)\neq 0.$ Thus, the Dynkin diagram of $\Psi_\aa(\SS)$ is obtained from the extended Dynkin diagram of $\SS$ by eliminating nodes representing the minimal root and all its neighbors.  Table \ref{1.4.I} follows from Figure \ref{A}.
\end{proof}

We deduce from \ref{1.4.I} by induction in $|O|$ that:
\smallskip
\item\label{1.4.J}
If  $O$ is an orthogonal subset of an irreducible root system $\SS$, then $\Psi_O(\SS)$ contains at most one irreducible component not of type $A_1$.

We denote this  irreducible component by~$\Theta_O(\SS)$ and we put ~$\Theta_O(\SS)=\emptyset$ if it does not exist.

\end{asm}
Importance of \ref{1.4.J} for the classification of subsystems is illustrated by the following example.
\begin{example}
\label{ex_A2m}
We claim that $I(L,S)=W(L,S)$ for every irreducible projective root system
$S$ and every its $\Pi$-system $L=A_{2m}$. Let $e$ be an end of $L$ and let $f\in I(L,S)$.  Since $W$ is transitive on $S$, we may assume $f(e)=e$. The set $\Psi_e(L)$ is a $\Pi$-system of type $A_{2m-2}$ in the projective root system $\Psi_e(S)$. Hence $\Psi_e(L)\subset\Theta_e(S)$. We deduce that  $f\in W(L,S)$ by induction in $m$ and Reduction Lemma.
\end{example}

\subsection{Core groups}
\begin{defn}
Let $M$ be a moset of $\SS$. Group $W_M$ is called a \emph{core group} of $\SS$.
\end{defn}

It follows from Proposition \ref{moset_conj} that core groups corresponding to all mosets $M$ are isomorphic. Therefore the type of $W_M$ is determined by the type of $\SS$ and we denote it $\Cor(\SS)$. We use notation $\nu(\SS)$ for the order of this group.

A core group $W_M$ can be identified with a subgroup  of  the permutation group $\SP(M)$. Indeed, there is a natural homomorphism of $W_M$ to  $\SP(M)$. It follows from Proposition \ref{humphreys} that $w=\Id$ if $w\in W$ acts trivially on $M$ because  none of roots is orthogonal to $M$.

\begin{asm}
\item\label{1.9.A}
Let $\SS=\cup_k \SS_k$ be a partition  of $\SS$ into  irreducible components. Then  $M_k=M\cap\SS_k$ is a moset in $\SS_k$, and $W_M$ is the direct product of $W_{M_k}$.

\smallskip

\begin{prop}
\label{core}
If $\LL\subset M$ and $w(\LL)\subset M$ for some $w\in W$, then there exists a $w'\in W_M$ such that $w'(\aa)=w(\aa)$ for all $\aa\in\LL$.
\end{prop}

\smallskip
\begin{proof}
The set $w(\LL)$ belongs to mosets $M$ and $w(M)$. If follows from \ref{1.8.A} that $M_1=M\setminus{w(\LL)}$ and $M_2=w(M)\setminus{w(\LL)}$ are mosets of the orthogonal complement $\td\SS$ to $w(\LL)$ in $\SS$. By Proposition \ref{moset_conj}, there is a $\td{w}\in W(\td\SS)\subset W$ such that $\td{w}(M_2)=M_1$. Since $\td{w}$ acts trivially on $w(\LL)$, $w'=\td{w}w$ preserves $M$.
\end{proof}

\smallskip
\item\label{1.9.AA}
Suppose $L$ is a $\Pi$-system in $S$, $M$ is a perfect moset in $L$ and $f\in I(L,S)$.
Let $O$ be a moset in $L$ such that $f(O)\subset M$. Then $f\in W(L,S)$ if $f|_M\in W_M(O,M)$.

\smallskip

\item\label{1.9.D}
If $\SS$ is irreducible, then $W_M$ acts transitively on $M$.

\smallskip
This follows from Proposition \ref{root_conj} and from Proposition \ref{core} applied to singletons.
\smallskip
\item\label{1.9.E}
The isotropy group $(W_M)_\aa$, where $\aa\in M$, is isomorphic to the core group of $\Psi_\aa(\SS)$.

\smallskip
This follows from Proposition \ref{humphreys}.

We conclude from  \ref{1.9.D} and \ref{1.9.E}, that:

\smallskip
\item\label{1.9.G}
$
\nu(\SS)=\nu(\Psi_\aa(\SS))\mu(\SS).
$

\smallskip
\item\label{1.9.H}
We have:
\begin{longtable}{|c|c|c|c|c|c|c|}
\caption{Orders of core groups} \label{ORD}\\
\hline
$\SS:$ & $A_n$ & $D_{2m+1}$ & $D_{2m}$ & $E_6$ & $E_7$ & $E_8$ \\
\hhline{|=|=|=|=|=|=|=|}
$\nu(\SS):$ & $\mu!$ & $2^m\cdot  m!$ & $2^{m-1}\cdot  m!$ & $\mu!$ & $168$ & $1344$ \\
\hline
\end{longtable}
where $\mu=\mu(\SS)$.

This follows from \ref{1.9.G} and Tables  \ref{tocl}, \ref{mu}.
\end{asm}

Suppose $\SS$ is irreducible. Since $W_M\subset\SP(M)$,  \ref{1.9.H} implies:
\begin{thm}
\label{thm_AnE6}
$W_M=\SP(M)$ if and only if $\SS=A_n$ or $E_6$.
\end{thm}

Core groups of $D_n$, $E_7$ and $E_8$ will be described in terms of enhanced Dynkin diagrams. For  $D_{2m}$ and $D_{2m+1}$ we represent  $M$ by a matrix
\begin{equation}
\label{lbl-m}
\begin{pmatrix}
1&3&\dots & 2m-1 \\
1'&3'&\dots & (2m-1)'
\end{pmatrix}.
\end{equation}
For $E_7$ and $E_8$ a  new labeling of $M$ is introduced on Figure \ref{new} where the labels are vectors of the space $\FB_2^3$ \ over   the field $\FB_2$ \   with  elements $0,1$.

\begin{figure}[htp]
\begin{gather*}
\genfrac{}{}{0pt}{}{\includegraphics{F.1}}{\DD(E_7)}\qquad
\genfrac{}{}{0pt}{}{\includegraphics{Js.1}}{\DD(E_8)}
\end{gather*}
\caption{$\DD(E_7)$ and $\DD(E_8)$}\label{new}
\end{figure}

\begin{thm}
\label{thm_Dn}
The core group $\Cor(D_{2m+1})$ can be represented as the group of permutations of columns in the matrix \eqref{lbl-m} combined with transpositions of the entries in some of columns. Permutations with  an even number of such transpositions form a subgroup representing $\Cor(D_{2m})$.
\end{thm}

\begin{proof}
Transpositions $\tau_i$ of  $2i-1$ and $(2i-1)'$ leaving fixed the rest of the entries generate a subgroup $T$ of $\SP(M)$ of order $2^m$. Transpositions $\ss_i$ of columns containing $2i-1$ and $2i+1$  generate a subgroup $S$ of $\SP(M)$ of order $m!$. The pair $S,T$ generates a subgroup  of order $2^m\cdot m!$ that is equal to 
$W_M$ by \ref{1.9.H}. To prove
Theorem \ref{thm_Dn} for $D_{2m+1}$ it is sufficient to demonstrate that $\tau_i$ and $\ss_i$ belong to 
$W_M$.
Similarly, to prove Theorem \ref{thm_Dn} for $D_{2m+1}$ it is sufficient to demonstrate that $\tau_i\tau_{i+1}$ and $\ss_i$ belong to $W_M$.

First, we consider the case of $D_4$. In this case every permutation of $M$ is induced by an automorphism of $\DD(D_4)$ which, by Proposition \ref{aut_enhanced},  is induced by an automorphism of $\SS$. Hence $W_M$ is a normal subgroup of $\SP(M)$. There is only one such subgroup of four elements, and it is generated by $\ss_1,\ss_2$ and $\tau_1\tau_2$. The case of any $D_{2m}, m\ge 2$, follows from $D_4$.

In order to prove the Theorem for $D_{2m+1}$, it suffices to show $\tau_1\in W_M$. But $\tau_1$ is induced by an automorphism of $\DD$, hence the statement follows from Proposition \ref{aut_enhanced}.
\end{proof}

In the case of  $E_7$,  only nonzero vectors of $\FB_2^3$ are used for labeling $M$.

\begin{thm}
\label{thm_E7}
$\Cor(E_7)$ is represented by $\GL_3(\FB_2)$.
\end{thm}

More precisely, we identify $\Cor(E_7)$ with a subgroup $\SP_{cor}$ of $\SP(M)$ and we demonstrate that $\SP_{cor}$ coincides with the group $\SP_{gl}=
\GL_3(\FB_2)$ of all invertible linear transformations of  $\FB_2^3$.

\begin{proof}
Note that $|\SP_{gl}|=168$.  By  \ref{1.9.H},  $\SP_{cor}$ has the same order.
Therefore it is sufficient to prove $|\SP|\ge 168$ for the group $\SP=\SP_{cor}\cap\SP_{gl}$.

Every automorphism $\ss$ of $\DD(E_7)$ preserves the vector $111$. We see from Figure \ref{new}  that  $x+y+z=0$  if and only if  the nodes with these labels lie  either on the same face of the tetrahedron or on a line through its center. This property is preserved under $\ss$ and therefore  $\ss(x)+\ss(y)+\ss(z)=0$. Since $z=x+y$, we have $\ss(x+y)=\ss(x)+\ss(y)$, which implies  $\ss\in\SP_{gl}$.
By Proposition \ref{aut_enhanced}, every automorphism $\ss$ is induced by an element of the Weyl group which implies that automorphisms $\ss_1,\dots,\ss_{24}$ determined by 24 rotations of the tetrahedron belong to $\SP_{cor}$.

Now we consider an automorphism of  $\hat{D_4}$-set $\{2,3,4,5,\ell_1\}\subset\DD(E_7)$ on Figure \ref{Fig} which transposes 3 with $\ell_1$ and 2 with 5. By Theorem \ref{thm_Dn} it is induced by a $\tau\in \SP_{cor}$  not moving elements $7,\ell_3,\ell_4$. We use Figure \ref{new} to check that $\tau(x+y)=\tau(x)+\tau(y)$ for all $x,y\in M$. Therefore $\tau\in\SP_{gl}$.

Since $|M|=7$ is prime, transformations $\tau,\ss_1,\dots,\ss_{24}$ generate a subgroup of $\SP$  acting transitively on $M$. The order of this subgroup is equal at least to $24\times 7=168$.
\end{proof}

\begin{defn}
\label{parity}
The parity $\rho(O)$ of a subset $O\subset M$ is equal to 0 if the sum of all labels of $O$ is \textbf{0} and equal to 1 if this sum is not \textbf{0}. 
\end{defn}
\begin{corol}
\label{orth_E7}
Two subsets of $M$ with the same cardinality $k$ are congugate if $k\neq 3,4$. If $k=3$ or $4$, then they are conjugate if and only if they have the same parity.
\end{corol}

It is sufficient to check this for $k=1,2,3$ because  $w(M\setminus{O})=M\setminus{O'}$ if $w(O)=O'$. Note that, for $|O|=3$, $\rho(O)=1$ if and only if $O$ corresponds to a basis of $\FB_2^3$.  The rest follows from Theorem \ref{thm_E7}. 

\begin{corol}
\label{aut_orth_E7}
 $W_O$ induces on $O$ the group of all permutations if $|O|=3$ or if $|O|=4$ and $\rho(O)=0$. If $|O|=4$ and $\rho(O)=1$, then  $W_O$ induces on $O$ the group of all permutations preserving an element $\aa$.\footnote{$\aa$ is a unique element such that $\rho(O\setminus\{\aa\})=0$.}
\end{corol}

\begin{figure}[htp]
\begin{gather*}
\genfrac{}{}{0pt}{}{\includegraphics{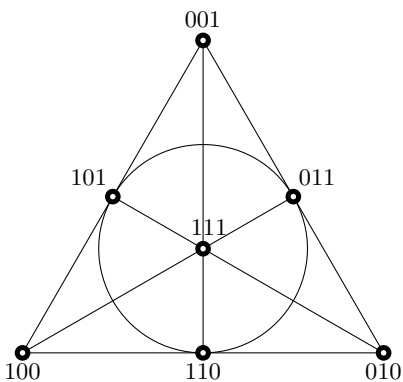}}{}
\end{gather*}
\caption{Fano plane}\label{Fano}
\end{figure}

\begin{rmrk}
By comparing Figure \ref{new} with Figure \ref{Fano} we see that  a moset $M$ can be interpreted as a projective plane over $\FB_2$ (Fano plane). Linearly dependent triplets are represented by lines in the plane and $\Cor(E_7)$ is represented by the  collinearity group.
\end{rmrk}

Mosets $M$ and $M^\diamond$ of $E_7$ and $E_8$ on Figure \ref{new} are related by the formula
\[
M^\diamond=M\cup\{\textbf{0}\}
\]
 where \textbf{0} is the zero vector 000.
\begin{thm}
\label{thm_E8}
$\Cor(E_8)$ is represented by the group $\AFF_3(\FB_2)$ of all affine transformations of $\FB_2^3$.
\end{thm}
\begin{proof}
Denote $\SP^\diamond_{aff}$ and $\SP_{cor}^\diamond$ subgroups of $\SP(M^\diamond)$ representing
$\AFF_3(\FB_2)$ and $\Cor(E_8)$. Group $\SP_{gl}=\SP_{cor}$ introduced in proof of Theorem \ref{thm_E7} is naturally embedded into $\SP^\diamond_{aff}\cap\SP_{cor}^\diamond$. By \ref{1.9.E} and Theorem \ref{thm_E7},  the isotropy group  of $\SP_{cor}^\diamond$ at \textbf{0} coincides with $\SP_{cor}$.
The parallel translation $\tau$ by 111 (in $M^\diamond$) is an element of $\SP^\diamond_{aff}$. It belongs also to $\SP_{cor}^\diamond$ because it
 is induced by an automorphism $\ff$ of $\DD(E_8)$.
\footnote{To define $\ff$ we relabel  $\DD(E_8)$ by elements $(a,b)$ of  group $\FB_4\times\FB_4$ and we put $\ff(a,b)= (a+2,b+2)$.}

Group $\SP^\diamond_{aff}$ is generated by  $\SP_{gl}$ and $\tau$.
Group generated by $\SP_{cor}$ and $\tau$ acts transitively on $M^\diamond$ and therefore its order is at least equal to $8|\SP_{cor}|=8\cdot 168=|\SP_{cor}^\diamond|$.
Since $\SP_{gl}=\SP_{cor}$, we conclude that $\SP^\diamond_{aff}\subset\SP^\diamond_{cor}$, and $\SP^\diamond_{aff}=\SP^\diamond_{cor}$ because both groups have the same order.
\end{proof}

\begin{corol}
\label{orth_E8}
Two subsets of $M$ with the same cardinality $k$ are congugate if $k\neq 4$. If $k=4$, then they are conjugate if and only if they have the same parity.
\end{corol}

Proof of Corollary \ref{orth_E8} is similar to that of Corollary \ref{orth_E7}.

A parity $\rho(O)$ of an orthogonal set $O$ is a parity of a $O'\subset M$ conjugate to $O$. (By Corollaries \ref{orth_E7} and  \ref{orth_E8}, $\rho(O)$ does not depend on the choice $O'$.) 

\begin{corol}
\label{aut_orth_E8}
If $|O|=4$, then $W_O$ induces on $O$ the group of all permutations.
\end{corol}

\begin{rmrk}
Consider a function $\eta$  defined on irreducible ADE root systems  by a table:

\begin{longtable}{|c|c|c|c|c|c|}
\hline
$\SS:$ & $A_n$ & $D_n$ & $E_6$ & $E_7$ & $E_8$\\
\hline
$\eta(\SS):$ & $\mu(A_n)$ & 1 & $\mu(E_6)$ & 2& 3\\
\hline
\end{longtable}

\textit{The condition $\eta(\SS)<|O|<\mu(\SS)-\eta(\SS)$ is necessary and sufficient for existence of a set $O'$   isomorphic but not conjugate to an orthogonal set $O$.}
\end{rmrk}
We do not use this fact.
 It can be proved by using Theorem \ref{thm_Dn}, Corollaries \ref{orth_E7}, \ref{orth_E8} and Table \ref{tocl}. [$\eta(\SS)$ is a minimal value of $|O|$ such that $\Psi_O(\SS)$ is empty or reducible.]

\subsection{Orthogonal subsets of an enhanced basis}
Suppose $S$ is a projective root system. We want to classify Weyl orbits in $\PP$ by using a classification of orbits of  a core group $W_M$. It is convenient to work not with the action of $W$ on $\PP$ but with its action on embeddings of $L\in\PP$ into $S$.
Maps $f_1,f_2\in I(L,S)$ lie on the same Weyl orbit if and only if $f_2f_1^{-1}\in W(f_1(L),S)$. For a moset $M\subset S$ and a subset $O\subset M$ put $f\in W_M(O,M)$ if $f=w$ on $O$ for a $w\in W_M$.
It  follows from Proposition \ref{core} and Theorem \ref{thm_perfect} that:

\begin{prop}
\label{reduction_orth}
Let $L\in\PP$ and $f\in I(L,S)$. For a perfect moset $O$ of $L$, a moset $M$ of $S$  and any embeddings $f_1\in W(O,M)$, $f_2\in W(f(O),M)$ we have
\[
f\in W(L,S)\quad\text{if and only if}\quad f_2ff_1^{-1}\in W_M(f_1(O),M).
\]
\end{prop}

Let $\Phi\subset S$ be a projective enhanced basis containing $M$. We give an algorithm to determine which $f\in I(L,\Phi)$ belong to $W(L,\Phi)$. By Proposition \ref{reduction_orth}, it is sufficient for every orthogonal subset $O\subset\Phi$ to construct an embedding $f\in W(O,M)$.

Let $O\subset\Phi$ be an orthogonal subset, $O_2=O\cap M$ and $O_1=O\setminus{O_2}$. The set $\Phi_1=\Psi_{O_2}(\Phi)$ contains $O_1$ and $M_1=M\cap\Phi_1$.

\begin{prop}
\label{complement}
There is a subsystem $S_1\subset S$ such that $\Phi_1$ is its enhanced basis and $M_1$ is its moset.
\end{prop}

\begin{proof}
For  a singleton $O_2$ this is clear from Figure \ref{Fig}. For any $O_2$, we apply induction in $|O_2|$.
\end{proof}

Proposition \ref{complement} reduces the case of any $O$ to  $O$ disjoint from  $M$. It is clear from  Figure \ref{Fig} that $\Phi\setminus{M}$ is an orthogonal set. Therefore it is sufficient to find $f$ for $O=\Phi\setminus{M}$. In Table \ref{tab} we use labels on Figure \ref{Fig} to describe $f(a)$ for every $a\in O$.

\begin{longtable}{|c|c|c|c|}
\caption{Embeddings $f$ for irreducible $S$} \label{tab}\\
\hline
$S:$ & $A_{2m}\,\,\text{or}\,\, A_{2m+1}$ & $D_{2m}$ & $D_{2m+1}$ \\
\hline
$a\in O:$ & $2,4,\dots,\,\,\,\,\,2m\,\,\,\,\,$ & $2,4,\dots,2m-2$ & $2,4,\dots,\,\,\,\,\,2m\,\,\,\,\,$ \\
\hline
$f(a):$ & $1,3,\dots,2m-1$ & $1,3,\dots,2m-3$ & $1,3,\dots,2m-1$ \\
\hhline{|=|=|=|=|}
$S:$ & $E_6$ & $E_7$ & $E_8$\\
\hline
$a\in O:$ &  $1,4,6,\ell_2$ & $1,4,6,\ell_2$ & $1,4,6,8,\ell_2,\ell_6,\ell_7,\ell_8$ \\
\hline
$f(a):$ &  $3,\ell_1,5,2$ & $3,\ell_1,5,2$ & $3,\ell_1,5,\ell_5,2,\ell_4,7,\ell_3$ \\
\hline
\end{longtable}

To get this table we use  Proposition \ref{aut_enhanced}. We note that, for
$A_{2m}, E_6, E_8$,  $f$ is induced by an automorphism of $\Phi$  which is induced by a $w\in W$. For the rest of irreducible $S$ we find an enhanced basis $\td\Phi\subset\Phi$ of a subsystem $\td{S}$ of $S$ such that $O\subset\td\Phi$ and $\td\Phi\setminus{O}$ is a moset of $\td{S}$. Values of  $\td\Phi$ are described in the following table:

\begin{longtable}{|c|c|c|c|c|}
\caption{Set $\td\Phi$} \label{tab5}\\
\hline
$S:$ & $A_{2m+1}$ & $D_{2m}$ & $D_{2m+1}$ & $E_7$ \\
\hline
$\mathstrut\td\Phi:$ & $\{1,2,\dots,2m\}$ & $\{1,2,\dots,2m-2\}$ & $\{1,2,\dots,2m\}$ & $\{1,2,3,4,5,6,\ell_1,\ell_2\}$\\
\hline
\end{longtable}

\begin{example}
We use diagram $\DD(E_8)$ on Figure \ref{Fig} (representing an enhanced basis $\Phi$) to define a $\Pi$-system $L$ of type $A_7$ and its embedding $f$ into $\Phi$ such that $f$ is in $I(L,\Phi)$ but not in $W(L,\Phi)$. Recall that boldfaced nodes represent a moset $M$ for $E_8$. We put
\[
L=\{2,4,5,6,7,8,\ell_5\},\quad f(L)=\{3,1,\ell_1,\ell_2,2,\ell_7,\ell_5\},\quad O=\{2,5,7,\ell_5\}.
\]
Since $O$ and $f(O)$ have different parity, $f\not\in W_M(O,M)$.  By Proposition \ref{reduction_orth} applied to  $f_1=f_2=Id$,  $f\not\in W(L,\Phi)$.
\end{example}

\begin{example}
Now we define $L\in\PP(E_7)$  of type  $A_3+A_1$  and its embedding $f$ into $\Phi$ of class $W(L,\Phi)$. We take
\[
L=\{7,6,\ell_3,4\},\quad f(L)=\{1,3,4,6\},\quad O=\{7,\ell_3,4\}.
\]
By applying Proposition \ref{complement} and Table \ref{tab}, we find $f_1\in W(O,M)$ with image $\{7,\ell_3,5\}$ and $f_2\in W(f(O),M)$ with image $\{3,\ell_1,5\}$. Since $f_1(O)$ and $f_2(f(O))$ have the same parity,  $f_2ff_1^{-1}\in W(f_1(O),M)$ by Theorem \ref{thm_E7}. Finally, Proposition \ref{reduction_orth} implies $f\in W(L,\Phi)$.
\end{example}

\begin{rmrk}
It can be shown that every projective $\Pi$-system in $S$ is conjugate to a $L\subset\Phi$ with $O\subset M$.
\end{rmrk}

\section{Classification theorems}
If $X,Y$ are root sets or projective root sets, then we write $X\simeq Y$ if  they are isomorphic  and $X\sim Y$ if  they are conjugate.

\begin{defn}
A $\Pi$-system $\LL\subset\SS$ is called \emph{normal}, if $\LL\simeq\LL'$ implies $\LL\sim\LL'$ for every $\LL'\in\PP$. Otherwise $\LL$ is called \emph{special}.
\end{defn}

\begin{defn}
The \emph{significant part} of a $\Pi$-system $\LL\subset\SS$ is the union of all its irreducible components of types $A_{2m-1}$. We call $\LL$ \emph{significant} if it coincides with its significant part.
\end{defn}

Note that every significant $\Pi$-system contains a unique perfect moset. It is insignificant if at least one of its irreducible component is not of type $A_{2m-1}$.

In this section we prove the following theorems.
\begin{thm}
\label{thm_main(1)}
Two $\Pi$-systems in an irreducible $\SS$ are conjugate if and only if their significant parts are conjugate.
\end{thm}

\begin{thm}
\label{thm_main(2)}
All special $\Pi$-systems in an irreducible $\SS$ are significant unless $\SS$ is of type $D_n$ with $n\ge 7$.
\end{thm}

\begin{thm}
\label{thm_main(3)}
Suppose an irreducible $\SS$ is not of type $D_n, n\ge 7$. Then two significant $\Pi$-systems in $\SS$ are conjugate if and only if their perfect mosets are conjugate.
\end{thm}
To be able to refer to Lemmas \ref{lem_red} and \ref{lem_red1}, we prove these theorems restated  in terms of  projective root systems $S$ and embeddings $f\in I(L,S)$.

\subsection{Proof of Theorem \ref{thm_main(1)}}
\begin{asm}
\item\label{significant1}
Suppose  an irreducible  insignificant $L\in\PP$ is not of type $D_{2m}$. Then
$I(L,S)=W(L,S)$.

\begin{proof}
For $L=A_{2m}$ this follows  from Example \ref{ex_A2m}.

If $L\neq A_{2m}$, then it is of type $D_{2m+1}$ with $m>1$ or $E_m$ with $m=6,7,8$. We use an induction in $|L|$. Let $L^e$ be obtained by eliminating from $L$ an end $e$. By excluding from consideration  one of three ends in each $L$,
\footnote{ In notation of  Figure \ref{A}, these are $2m$ for $D_{2m+1}$,  1 for $E_7$ and  2 for $E_6$ and $E_8$.}
we obtain an insignificant $L^e\neq D_{2m}$, and, by induction assumption,  $f'=wf\vert_{L^e}=\Id$ for some  $w\in W$. If $\<f'(e)|e\>\neq 0$, then the diagram of $L^e\cup\{e,f'(e)\}$ is acyclic and, by  Proposition \ref{section}, it must be on Figure \ref{A}. With our selection of $e$, this is not the case. Hence $\<f'(e)|e\>\neq 0$ and, by Lemma \ref{lem_red1}, $f'$ and $f$ are in $W(L,S)$.
\end{proof}

\item\label{significant2}
If $L\in\PP$ is of type $D_{2m}$ and if $f\in I(L,S)$, then $fs\in W(L,S)$ for an automorphism $s$ of $L$.

\begin{proof}
Let $e_1$ and $e_2$ be ends of $L$ such that $\td{L}=L\setminus\{e_1,e_2\}$ is of type $A_{2m-2}$.  By \ref{significant1}, we may assume $f\vert_{\td{L}}=\Id$.
The set $E=\{e_1,e_2,f(e_1),f(e_2)\}$ is not orthogonal. Indeed, otherwise an acyclic diagram of $\td{L}\cup E$ would have 5 ends, which contradicts Proposition \ref{section}. Without loss of generality, we assume $\<e_2|f(e_2)\>\neq 0$ or $\<e_2|f(e_1)\>\neq 0$. In any case, there is a $q\in\Aut{L}$, trivial on $\td{L}$, such that $\<e_2|fq(e_2)\>\neq 0$. By Lemma \ref{lem_red1}, $wfq\vert_{L_1}=\Id$ for $L_1=L\setminus\{e_1\}$ and a $w\in W$.

Set $f'=wfq$. If $\<f'(e_1)|e_1\>\neq 0$, then $f'\in W(L,S)$ by Lemma \ref{lem_red1}. Hence $fs\in W(L,S)$ for $s=q$. Case $\<f'(e_1)|e_1\>=0$ is possible only if $m=2$. Then $\hat{L}=L\cup\{f'(e_1)\}$ is the extended projective $\Pi$-system and there is a $w'\in W_{\hat{L}}$ transposing $e_1$ and $f'(e_1)$ and preserving $L_1$. Hence $w'f'(L)=L$ and $w'f'$ induces on $L$ an  automorphism $q_1$, and \ref{significant2} holds for $s=qq_1^{-1}$.
\end{proof}

\end{asm}
\begin{prop}
\label{significant}
Suppose $A$ is the significant part of  $L\in\PP$ and $f(A)\sim A$. Then $fs\in W(L,S)$ for an automorphism $s$ of $L$ preserving each irreducible insignificant component   and trivial on all of them not of type $D_{2m}$.
\end{prop}
\begin{proof}
We use induction in the number of irreducible insignificant components of $L$. Let $Q$ be one of them. Put $L'=L\setminus{Q}$. Suppose $f\vert_{L'}=\Id$. Let $O'$ be a perfect moset of $L'$. Then $Q$ and $f(Q)$ are irreducible insignificant  $\Pi$-systems in irreducible projective root system $S'=\Theta_{O'}(S)$. There is an automorphism $q$ of $Q$ such that $w'f'q=\Id$ where $f'=f\vert_{Q}$ and $w'\in W$. Extend $q$ trivially to an automorphism $s$ of $L$, and set $f''=w'fs$. If $O\supset O'$ is a perfect moset of $L$, we have $f''\vert_{O}=\Id$. By Theorem \ref{thm_perfect}, $f''\in W(L,S)$, hence $fs\in W(L,S)$.
\end{proof}
\subsection{Proof of Theorem \ref{thm_main(2)}}
We deal with the class $\SF$ of all irreducible projective ADE systems $S$ except $D_n, n\ge 7$. Set $\AF$ to be the class of projective ADE root systems having all special $\Pi$-systems significant. In order to prove Theorem \ref{thm_main(2)}, we have to show $\SF\subset\AF$. We  use the following    propositions.

\begin{prop}
\label{aux}
Suppose $L\subset S$ is a significant $\Pi$-system with the perfect moset $O$ and $f\in I(L,S)$. Let $L_1\subset L$ be an irreducible component with the perfect moset $O_1$ and $L_2=L\setminus{L_1}$. The following statements are equivalent:
\begin{enumerate}
 \item[(1)] there is a $w\in W$ such that $wf(L)=L$ and $wf(L_1)=L_1$;
 \item[(2)] there is a $w'\in W$ such that $w'f(L_2)=L_2$ and $w'f(O_1)=O_1$.
\end{enumerate}
\end{prop}

\begin{proof}
 $(1)\Rightarrow(2)$ with $w'=w$. Suppose that (2) holds. By substituting  $w'f$ for $f$ we can assume that $f(L_2)=L_2$ and $f(O_1)=O_1$. Denote by $O_2$ the perfect moset of $L_2$. Note $S_1=\Psi_{O_2}(S)$ contains $L_1$ and $f(L_1)$. (Indeed, $f(L_1)\perp f(L_2)=L_2\supset O_2$.) Since $L_1=A_k$, by Theorem \ref{thm_AnE6}, there is a $w_1\in W(S_1)$ such that $w_1f=\Id$ on $O_1$ (and $w_1=\Id$ on $O_2$). Let $h$ be an automorphism of $L$ such that $h=\Id$ on $L_1$ and $fh=\Id$ on $L_2$. Then we have $w_1fh=\Id$ on $O_1$ and $O_2$, and therefore $w_1fh=\Id$ on $O=O_1\cup O_2$. By Theorem \ref{thm_perfect}, $w_1fh=w_2$ on $L$, where $w_2\in W$. We put $w=w_2^{-1}w_1$. Note that $wfh=\Id$ on$L$ and $w=\Id=hf$ on  $L_2$. Therefore  $wf(L)=wfh(L)=L$ and $wf(L_1)=wfh(L_1)=L_1$.
\end{proof}

\begin{prop}
\label{transf_psi}
If $S$ is irreducible and $\Psi_a(S)\in\AF$ for all $a\in S$, then $S\in\AF$.
\end{prop}

\begin{proof}
Assuming that $\Psi_a(S)\in\AF$ for all $a\in S$ we check that, if $L\in\PP(S)$ has the significant part $A\neq L$, then $f(L)\sim L$ for all  $f\in I(L,S)$. This follows from Theorem \ref{thm_main(1)} if $A=\emptyset$.

If $A$ is non-trivial, then we consider its irreducible component $A_1$  of minimal rank and an end $e$ in $A_1$. Replacing $f$ by $wf$, $w\in W$, we can assume $f(e)=e$. $\Pi$-system $\tilde{L}=\Psi_e(L)$ in $\tilde{S}=\Psi_e(S)$ is insignificant. Hence, by our assumption, $f(\tilde{L})=\td{w}(\tilde L)$ for a $\td{w}\in W(\tilde{S})$. The irreducible component $\td{A_1}=\Psi_e(A_1)$ of $\td{L}$ has the smallest rank, and therefore is preserved by $\td{w}$. Then $\td{w}$ preserves $A\setminus{A_1}$ and the perfect moset of $A_1$. By Proposition \ref{aux}, $f(L)\sim L$.
\end{proof}

Denote by $\SF_1$ the class $\SF$  augmented by reducible types $A_2+A_1$, $A_3+A_1$ and $D_4+A_1$ and by all orthogonal types $kA_1, k\ge 1$.

\begin{asm}
\item\label{5.3.A}
If $S\in\SF_1$, then $\Psi_a(S)\in\SF_1$ for all $a\in S$.
This  follows immediately from Table~\ref{tocl}.

\smallskip
\item\label{5.3.B}
Any reducible $S\in \SF_1$ belongs to $\AF$.
\begin{proof}
The statement is trivial for orthogonal~$S$. Suppose  $S=A_2+A_1$, $A_3+A_1$ or $D_4+A_1$ and let $L\in\PP(S)$ be insignificant. If $L$ is irreducible, then $L$ must belong to a unique irreducible component of $S$ not of type $A_1$, and it is normal by Theorem \ref{thm_main(1)}. If $L$ is reducible, then $L\supset\tilde{L}$ where  $\tilde{L}$ is of type $A_2+A_1$. [$A_2$ is a subset of an irreducible component of $L$ and $A_1$ is contained in another component of $L$.] But none of $A_2$, $A_3$, $D_4$ contains $A_2+A_1$, hence the singleton in $\tilde{L}$ coincides with that in~$S$. On the other hand, $A_2$ is a maximal proper insignificant $\Pi$-system in $A_3$ and~$D_4$. Hence $L=\tilde{L}$ and the statement follows from normality of $A_2$.
\end{proof}
\end{asm}

\begin{proof}[Proof of Theorem \ref{thm_main(2)}]
We will prove that $\SF_1\subset\AF$.\footnote{In fact, $\SF_1=\AF$. We do not use this fact.} We use induction in rank of $S\in\SF_1$. In case of reducible $S$ we have $S\in\AF$ by~\ref{5.3.B}. If $S$ is irreducible, we use Proposition~\ref{transf_psi} and \ref{5.3.A}.
\end{proof}

\subsection{Proof of Theorem \ref{thm_main(3)}}
Theorem \ref{thm_main(3)} is equivalent to the statement:

\begin{prop}
\label{rigid}
Let $S\in\SF$, $L\in\PP(S)$ and $f\in I(L,S)$. If $f(O)=O$ for a perfect moset $O$ of $L$, then $f(L)\sim L$.
\end{prop}

Proposition \ref{rigid} is trivial if $L$ is normal and it follows from Theorem \ref{thm_perfect} if $W_O=\SP(O)$. By Theorem \ref{thm_AnE6}, it remains to prove it for a special $L$ and $S$ of type $D_4$, $D_5$, $D_6$, $E_7$, $E_8$. By Theorem \ref{thm_main(2)}, $L$ is significant. In the proof we use two lemmas.

\begin{lem}
\label{lem_a3a3}
If $L\in\PP(S)$ is of type $2A_3$, then there exists a $w\in W$ transposing the irreducible components of $L$.
\end{lem}

\begin{proof}
Only $D_6$, $E_7$ and $E_8$ contain such $L$. We have $|O|=4$. Let $G_O$ be the group of transformations of $O$ induced by $W_O$. It follows from Theorem \ref{thm_perfect} that if $G_O=\SP(O)\simeq\SP_4$, then $I(L,S)=W(L,S)$. For $S=E_8$, by Corollary \ref{aut_orth_E8}, $G_O\simeq\SP_4$. Corollary \ref{aut_orth_E7} implies that, in the case of $E_7$, $G_O$ is isomorphic either to $\SP_4$ or to $\SP_3$. However we deduce from Proposition \ref{aut_enhanced} (applied to $\Phi$ of type $2A_3$) that $G_O$ contains a subgroup isomorphic to $\FB_2\times\FB_2$. Hence $G_O\simeq\SP_4$.

A unique  subdiagram of type $2A_3$ of  diagram $\DD(D_6)$ is
$\{1,2,1',5,4,5'\}$ (in notation of Figure \ref{Fig}). Hence, $O=\{1,1',5,5'\}$. Let $f$ be an automorphism of $L$ inducing on $O$ transpositions of $1$ with $5$ and of $1'$ with $5'$. By Theorem \ref{thm_Dn}, $f\in W(L,S)$.
\end{proof}

Denote by $\SF_2$ the class of projective root systems $S=S_0+S_1$, where $S_0\in\SF$ and $S_1$ is of type $kA_1$. Note that $\Psi_a(S)\in\SF_2$ if $S\in\SF_2$, for all $a\in S$.

\begin{lem}
\label{lem_A3kA1}
Let $S\in\SF_2$, $L\in\PP(S)$ have type $A_3+kA_1$, $f\in I(L,S)$. If $f(O)=O$ for a perfect moset $O$ of $L$, then $f(L)\sim L$.
\end{lem}
\begin{proof}
Denote by $L_1\subset L$ the irreducible component of type $A_3$. By Proposition \ref{aux}, it suffices to find a $w\in W$ such that $wf(O)=O$ and $wf(O_1)=O_1$ where $O_1$ is the perfect moset of $L_1$.

Put
\[
Q=L_1\cup f(L_1),\quad I=Q\cap O,\quad J=O\setminus I.
\]
Note that $Q\perp J$. Since $Q$ does not contain any singleton, it belongs to a unique irreducible component $\td{S}$ of $\Psi_J(S)$ (see \ref{1.4.J}).

Put  $K=O_1\cap f(O_1)$. If $|K|=0$, then $|O_1\cup f(O_1)|=4$, $|Q|=5$ or $6$ and $Q$ is an acyclic set of type $2A_3$, $\hat{D_4}$ or $\hat{D_5}$. In the first case, by Lemma \ref{lem_a3a3}, there is a $w\in W(\tilde{S})$ that transposes $L_1$ and $f(L_1)$. We have $w(O)=O$, $wf(L_1)=L_1$ and $wf(O_1)=O_1$. For the rest of cases $O$ is a moset of a subsystem of type $D_4$ or $D_5$, and Theorem \ref{thm_Dn} implies $wf(O_1)=O_1$ for some $w\in W_O$.

If  $|K|=2$, then $f(O_1)=O_1$ and we take $w=\Id$.

Finally, if $|K|=1$, then $f(O_1)\cap O_1$ is a singleton $\{a\}$. By Theorem \ref{thm_AnE6}, applied to $A_3$, we can assume that $f(a)=a$. Let $O_1=\{a,b\}$. Then $I=\{a,b,c\}$ where $c=f(b)$. We consider $w_1,w_2\in W(\td{S})_I$ transposing respectively $a$ with $b$ and $a$ with $c$. Hence $wf(O)=O$ and $wf(O_1)=O_1$, where $w=w_1w_2$.
\end{proof}
		
\begin{proof}[Proof of Proposition~\ref{rigid}]
By induction in $|L|$, we prove  a more general statement: if  $S\in\SF_2$, then
for every irreducible component $L_1\subset L$ not of type $A_1$ there is a $w_1\in W$ such that $w_1f(L)=L$ and $w_1f(L_1)=L_1$.

This is trivial for  $L=kA_1$. Case  $L=A_3+kA_1$  follows from Lemma \ref{lem_A3kA1}. For the rest of types we consider an irreducible component  $L_1\neq A_1$ and we choose its end  $e$ and  take $L_1'\subset L_1$ of type $A_3$  that contains $e$. Since $L'=L_1'\cup O$ is of type $A_3+kA_1$, we conclude from  Lemma \ref{lem_A3kA1} that $wf(L')=L'$ for some $w\in W$. Hence $wf(L_1')=L_1'$ and $wf(O)=O$. By Theorem \ref{thm_AnE6}, applied to $A_3$, we can choose $w$ such that $wf(e)=e$.

Put $\td{S}=\Psi_e(S)$, $\td{L}=\Psi_e(L)$, $\td{L_1}=\Psi_e(L_1)$. We have $\td{L}, wf(\td{L})\subset\td{S}\in\SF_2$. If $\td{L_1}\neq A_1$, then by induction, there is a $w_1\in W(\td{S})$ such that $w_1wf(\td{L})=\td{L}$ and $w_1wf(\td{L_1})=\td{L_1}$. In particular, $w_1wf$ preserves $L\setminus{L_1}$ and the perfect moset of $L_1$. By Proposition \ref{aux},  $w_2w_1wf(L)=L$ and $w_2w_1wf(L_1)=L_1$ for a $w_2\in W$.

Finally, suppose $\td{L_1}=A_1$. Then $L_1=L_1'$ has type $A_3$ and $wf(L_1)=L_1$. Set $O_1$ to be the perfect moset of $L_1$, $L_2=L\setminus{L_1}$ and $S_1=\Psi_{O_1}(S)\in\SF_2$. Since $wf$ preserves the perfect moset of $L_2$, we obtain by induction that $w_1wf(L_2)=L_2$ for a $w_1\in W(S_1)$. Since $w_1wf(O_1)=O_1$, Proposition \ref{aux} implies the existence of a $w_2\in W$ such that $w_2w_1wf(L)=L$ and $w_2w_1wf(L_1)=L_1$.
\end{proof}

\subsection{Description of Weyl orbits}
To classify Weyl orbits in $\PP(S)$ we consider subdiagrams $L$ of $\DD(S)$ representing elements of $\PP(S)$ and we investigate which of them are normal and which of isomorphic special $\Pi$-diagrams are conjugate. The general case can be easily reduced to the case of irreducible $S$. For $S=A_n$ and $E_6$ a solution is given by:

\begin{thm}
\label{thm_normal}
All $\Pi$-systems in $A_n$ and $E_6$ are normal.
\end{thm}

This follows from Theorems \ref{thm_perfect} and \ref{thm_AnE6}.

\subsubsection{Case of $D_n$}
In the case of $D_n$, we call a $\Pi$-diagram\footnote{By a $\Pi$-diagram we mean a Dynkin subdiagram of $\DD(\SS)$.} $L\subset\DD(D_n)$ \textit{thin}\index{thin $\Pi$-diagram} if it contains no pair $\{i,i'\}$, and we call it \textit{thick}\index{thick $\Pi$-diagram} if it is not thin. Thick diagrams of type $2A_1$ are conjugate and the same is true for type~$A_3$. We call them $D_2$- and $D_3$-diagrams. The number of $D_i$-subdiagrams of $L$ is denoted by $\dd_i=\dd_i(L)$. We call $(\dd_2,\dd_3)$ the \textit{tag of}\index{tag of $\Pi$-diagram}~$L$. We denote by $\oo(L)$ and we call the \textit{width}\index{width of $\Pi$-diagram} of a diagram $L$ the minimal value of $k$ such that $L$ is contained in $\DD(D_k)$.
Thin significant diagrams of   width $n$ are called \textit{distinguished}\index{distinguished diagram}. Such diagrams exist only if $n$ is even.

\begin{thm}
\label{thm orb-Dn}
Every class of isomorphic distinguished diagrams consists of two W-orbits transformed into each other by any automorphism $\tau_i$ of $\DD(D_n)$ transposing $i$ and $i'$ and not moving the rest of nodes. In all other cases isomorphic $\Pi$-diagrams are conjugate if and only if they have the same tag.
\end{thm}

Since subdiagrams of any $L\subset\PP(D_n)$ of classes $D_2$  and $D_3$ are contained in the significant part $A$ of $L$, the tags of $L$ and $A$ coincide. Besides, all distinguished diagrams are significant. By Theorem \ref{thm_main(1)}, it is sufficient to prove Theorem \ref{thm orb-Dn} for significant $L$.

\begin{asm}
\item\label{5.2.A}
 Let $N=N_1\cup N_2$ be the partition of the set $N=\{1,2,3,\dots,2m-1\}$ into odd and even numbers.
Every distinguished subdiagram of $\DD(D_{2m})$ corresponds to a pair $(I,J)$, where $I\subset N_1$ and $J\subset N_2$, by the formula
\[
L_I^J=\tau_I(N\setminus J),
\]
where $\tau_I=\prod_{i\in I}\tau_i$.  Diagrams $L_{I_1}^{J_1}$ and $L_{I_2}^{J_2}$ are isomorphic if  $J_1=J_2$, and they are conjugate if and only if, in addition,  $|I_1|$ and $|I_2|$ are both odd or both even.

\smallskip
This follows from Theorem \ref{thm_Dn}.

\item\label{5.2.D}
Not distinguished significant $\Pi$-diagrams are conjugate if and only if they have the same tag.

\begin{proof}
 Suppose $L'=w(L)$ where $w\in W$. Then $L$ and $L'$ have the same tag because $W$ preserves the types of connected  components of a $\Pi$-diagram and  it preserves the classes of  thin and thick diagrams.

On the other hand, if isomorphic $L$ and $L'$ have the same tag, then it is possible to establish a 1-1 correspondence between their components of each class $D_2, D_3$ and between the  thin connected  components of each type $A_{2i-1}$. We get this way an isomorphism $f:L\to L'$. We can choose $f$  such that the restriction of $f$ to the perfect moset of $L$ coincides with a permutation of columns in the matrix \eqref{lbl-m}. This permutation  belongs to $W$ by Theorem~\ref{thm_Dn}. By Theorem \ref{thm_perfect}, $f\in W$ and $L\sim f(L)=L'$.
\end{proof}
\end{asm}

\subsubsection{Case of $E_7$}
We define the parity of $L\in\PP(S)$ as the parity of its perfect moset $O$. Note that the parity does not depend on the choice of  $L$ within one Weyl orbit. The same is true for $|O|$ which we call the \textit{charge} of $L$.

\begin{thm}
\label{thm 5.5}
In the case of $E_7$ every special orbit belongs to one of three types
\begin{equation}
\label{E7.1}
A_5,  A_3+A_1,  3A_1
\end{equation}
with charge 3 or to one of types
\begin{equation}
\label{E7.1a}
A_5+A_1,  A_3+2A_1,  4A_1
\end{equation}
with charge 4. Each of these types consists of two Weyl orbits: one with parity 0 and the other with parity~1.
\end{thm}

\begin{rmrk}
Both orbits of charge 3 can be represented by distinguished subdiagrams of any $D_6$-subdiagram of $\DD(E_7)$ (for instance, of $L=\{2,3,4,5,6,7\}$). A similar representation, with $D_6+A_1$ (for instance,  $L+\{\ell_4\}$) substituted for  $L$, is possible for orbits of charge~4.
\end{rmrk}

By Theorem \ref{thm_main(2)} every special $L$ is significant and therefore it has a unique moset~$O$. By Theorem \ref{thm_main(3)} and by Proposition \ref{core},  Weyl orbits are in a 1-1 correspondence with the core orbits. We conclude from Corollary \ref{orth_E7} that all special $\Pi$-systems $L$ have charge 3 or 4  and that the class of all isomorphic special $L$ consists of two $W$-orbits: one with parity 0, the other with parity 1.  All significant types with charge 3 are $3A_1$, $A_3+A_1$, $A_5$.

All significant types with charge 4 are $4A_1$, $A_3+2A_1$, $A_5+A_1$, $2A_3$, $A_7$. But the latter two types are normal. Indeed, in these cases the group of transformations of $O$ induced by $W_O$ contains a subgroup isomorphic to $\FB_2\times\FB_2$. Then by Corollary \ref{aut_orth_E7} the parity of $L$ must be $0$ and therefore $L$ is normal.
We arrive to the  list of six types presented in Table~\ref{SpecE_7}.

\begin{longtable}{|c|c|c|}
\caption{Special Weyl orbits for $E_7$}
\label{SpecE_7}\\
\hline
Type & Parity 0 & Parity 1 \\
\hhline{|=|=|=|}
$3A_1$ & $\{2\}+\{5\}+\{7\}$ & $\{3\}+\{5\}+\{7\}$\\
\hline
$A_3+A_1$ & $\{5,6,7\}+\{2\}$ & $\{5,6,7\}+\{3\}$\\
\hline
$A_5$	& $\{2,4,5,6,7\}$ & $\{3,4,5,6,7\}$\\
\hline
$4A_1$ & $\{3\}+\{5\}+\{7\}+\{\ell_4\}$ & $\{2\}+\{5\}+\{7\}+\{\ell_4\}$\\
\hline
$A_3+2A_1$ & $\{5,6,7\}+\{3\}+\{\ell_4\}$ & $\{5,6,7\}+\{2\}+\{\ell_4\}$\\
\hline
$A_5+A_1$ & $\{3,4,5,6,7\}+\{\ell_4\}$ & $\{2,4,5,6,7\}+\{\ell_4\}$ \\
\hline
\end{longtable}

\subsubsection{Case of $E_8$}
\begin{thm}
\label{thm 5.6}
All special $L\in\PP(E_8)$ belong to one of five types
\begin{equation}
\label{E8.1}
A_7, A_5+A_1, 2A_3, A_3+2A_1, 4A_1
\end{equation}
with charge 4. To each of these types there correspond two Weyl orbits: one with parity 0 and the other with parity~1.\footnote{Both orbits  can be represented by distinguished  subdiagrams of any $D_8$-subdiagram of $\DD(E_8)$ (for instance, of $L=\{2,3,4,5,6,7,8,\ell_5\}$).}
\end{thm}

Theorems \ref{thm_main(2)}, \ref{thm_main(3)} and Corollary \ref{orth_E8} imply that all special $L$ are significant and every special type is represented by two Weyl orbits with charge 4 and different parities.

Every significant $L$ with charge 4 belongs to one of types $4A_1$, $A_3+2A_1$, $2A_3$, $A_5+A_1$ and $A_7$. These types are special which follows from Table \ref{SpecE_8}.

\begin{longtable}{|c|c|c|}
\caption{Special Weyl orbits for $E_8$}
\label{SpecE_8}\\
\hline
Type & Parity 0 & Parity 1 \\
\hhline{|=|=|=|}
$4A_1$ & $\{2\}+\{5\}+\{7\}+\{\ell_5\}$ & $\{3\}+\{5\}+\{7\}+\{\ell_5\}$\\
\hline
$A_3+2A_1$ & $\{7,8,\ell_5\}+\{5\}+\{2\}$ & $\{7,8,\ell_5\}+\{5\}+\{3\}$\\
\hline
$2A_3$ & $\{7,8,\ell_5\}+\{2,4,5\}$ & $\{7,8,\ell_5\}+\{3,4,5\}$\\
\hline
$A_5+A_1$ & $\{5,6,7,8,\ell_5\}+\{2\}$ & $\{5,6,7,8,\ell_5\}+\{3\}$ \\
\hline
$A_7$ & $\{2,4,5,6,7,8,\ell_5\}$ & $\{3,4,5,6,7,8,\ell_5\}$ \\
\hline
\end{longtable}

\subsection{Order between $W$-orbits}
If $\SS$ is the root system of an algebra $\ga$, then a partial order between conjugacy classes of subalgebras  defined in the Introduction can be investigated in terms of  an order  between W-orbits in $\PP(\SS)$. Theorems \ref{thm orb-Dn}, \ref{thm 5.5}, \ref{thm 5.6} allow to determine if $O_1\prec O_2$ for any pair  of orbits.

Every normal orbit is determined by its type. For special orbits we use additional labels. For instance, we label twelve special orbits in $\PP(E_7)$  represented on Table 7  by a combination of their types and parities   (indicated by  a superscript 0 or 1 at a type symbol). E. g., writing $[A_5+A_1]^0$ means type $A_5+A_1$ and parity 0. 
 Partial order between the corresponding orbits [which  follows immediately from Table 7]   is demonstrated by   the following graph where an arrow is  directed from label of $O_2$ to label of $O_1$ if  $O_1\prec O_2$.
\begin{center}
\textbf{Order between special orbits for $E_7$}
\end{center}
$$
  \begin{CD}
  [A_5+A_1]^0  @>>> [A_3+2A_1]^0  @>>> [4A_1]^0 \\
  @VVV @VVV @VVV\\
  [A_5]^1      @>>> [A_3+A_1]^1   @>>> [3A_1]^1 \\
  @. @AAA @AAA\\
  [A_5+A_1]^1  @>>> [A_3+2A_1]^1  @>>> [4A_1]^1 \\
  @VVV @VVV @VVV\\
  [A_5]^0      @>>> [A_3+A_1]^0   @>>> [3A_1]^0 \\
  \end{CD}
$$
 A similar directed graph for $E_8$ consists of two isomorphic connected components. Labels in each component are of the same parity.

\begin{center}
\textbf{Order between special orbits for $E_8$} ($i=\,0$ or $\,1$)
\end{center}

$$
  \begin{CD}
  [A_7]^i    @>>> [2A_3]^i     @. \\
  @VVV @VVV \\
  [A_5+A_1]^i @>>> [A_3+2A_1]^i @>>> [4A_1]^i \\
  \end{CD}
$$
Analogous graphs can be easily drawn  for $D_n$ with not too big $n$. 

For any ADE system and normal $O_1, O_2$, a relation $O_1\prec O_2$ holds
if and only if  the Dynkin diagram of  $O_1$ is a subdiagram of the  enhanced Dynkin diagram of  $O_2$. This criterion works also  in the case of a special $O_2$ and normal $O_1$.  

It remains to consider the case of special $O_1$ and normal $O_2$. We illustrate an approach to this case on the following example. Suppose $\SS=E_8$ and labels of $O_1$ and $O_2$ are  $[4A_1]^0$ and 
$E_6$. Subdiagram $\{2,3,\ell_1,\ell_2\}$ of  $\DD(E_6)$ has parity 0  as a subdiagram of $\DD(E_8)$. Therefore  $O_1\prec O_2$. On the other hand,  $O'_1$ with  label $[4A_1]^1$  is not comparable with $O_2$ because $4A_1$ is normal in $E_6$ and therefore its  parity in $E_8$ is 0.


\appendix
\section{Appendix}
\subsection{Basic facts}
A finite set $\SS$ of nonzero vectors in a Euclidean space $V$ is called a \emph{root system} if:
\begin{enumerate}
\item[(1)]
\label{rs1}
$\SS$ spans $V$;
\item[(2)]
\label{rs2}
for every $\aa,\bb\in\SS$,
$$
(\aa|\bb)=\frac{2(\aa,\bb)}{(\bb,\bb)}
$$
is an integer;
\item[(3)]
\label{rs3}
for every $\aa,\bb\in\SS$,
$$
s_\bb(\aa)=\aa-(\aa|\bb)\bb\in\SS;
$$
\item[(4)]
\label{rs4}
for every $\aa\in\SS$, $k\aa\not\in\SS$ for all $k\neq\pm 1$.
\end{enumerate}

Condition (4) not always is a part of the definition of root systems (see, e.g.,\cite{Ser01}). What we call a root system, many authors call a reduced root system.

The \emph{Weyl group} $W$ of $\SS$ is a subgroup of $\SP(\SS)$ generated by all $s_\aa$, $\aa\in\SS$. An \emph{automorphism} of $X\subset\SS$ is an element of $\SP(X)$ preserving the inner product $(\cdot,\cdot)$. We have $W\subset\Aut(\SS).$ A \emph{basis} of $\SS$ is a subset $\Pi\subset\SS$ such that every $\aa\in\SS$ has a unique presentation $\sum_{\gg\in\Pi}k_\gg\gg$, where all $k_\gg$ are non-negative integers or all non-positive.  Any root system contains a basis.

 Let $\SS$ be a root system and $W$ its Weyl group. We say that $X,Y\subset\SS$ are conjugate, if $w(X)=Y$ for some $w\in W$.

\begin{prop}
\label{base_conj}
All bases of $\SS$ are conjugate (\cite[Theorem 1.4]{Hum90}).
\end{prop}

\begin{prop}
\label{root_conj}
If $\SS$ is irreducible, then all roots of $\SS$ having the same length are conjugate (\cite[Chapter~VI, Proposition~11]{Bou02}).
\end{prop}

\begin{prop}
\label{humphreys}
If $w(\aa)=\aa$, where $\aa\in\SS$ and $w\in W$, then there are roots $\aa_1,\aa_2,\dots,\aa_k$ orthogonal to $\aa$, such that $w=s_{\aa_1}s_{\aa_2}\cdots s_{\aa_k}$ (\cite[Theorem 1.12]{Hum90}).
\end{prop}

\begin{prop}
\label{automorphisms}
Let $\Pi\subset\SS$ be a basis. All automorphisms of $\Pi$ extend uniquely to automorphisms of $\SS$, forming a subgroup $A\subset\Aut\SS$. The group $\Aut(\SS)$ is a semidirect product of $A$ and the normal subgroup $W$ (\cite[Theorem 1]{Dyn51}).
\end{prop}

\begin{prop}
\label{involution}
If $\SS$ is irreducible not of type $D_{2n}$, $n\ge 2$, then $\Aut\SS$ is generated by $W$ and $-\Id$ (\cite[Lemma 5]{Dyn51} or \cite[Theorem 0.16]{Dyn52a}).
\end{prop}

\begin{prop}
\label{Phi}
Suppose $\LL$ is an extended basis of an irreducible $\SS$.  If $\LL\setminus\{\aa\}$ and $\LL\setminus\{\bb\}$ are bases, then $w(\LL)=\LL$ and $w(\aa)=\bb$ for some 
$w\in W$.
\end{prop}

This follows from the list of  automorphisms for extended Dynkin diagrams  and Proposition \ref{automorphisms}.

\subsection{Notation Index}
\begin{flushleft}
$\Aut(X)$:  group of automorphisms of $X$\\

$f\vert_X$: restriction of a map $f$ to $X$\\
$G_X$:  subgroup  $\{g\in G:g(X)=X\}$ of a group $G$\\

$\Id$:  identity map\\
$I(X,\SS)$:  set of all embeddings of $X$ into $\SS$\\

$p(\aa)$:  projective root $(\aa,-\aa)$ corresponding to  a root $\aa$\\
$\PP$:  set of all $\Pi$-systems (in $\SS$ or in $S$)\\

$S:$  projective root system \\
$\SP(X)$:  group of all permutations on a finite set $X$\\
$W:$  Weyl group\\
$W(X,\SS)$:  set of all $f\in I(X,\SS)$ such that $f|_X=w|_X$ for some $w\in W$\\

$\GG_X$:  diagram of $X$\\
$\DD_X$:  projective diagram of $X$\\

$\DD(\SS)$: enhanced Dynkin diagram of $\SS$\\

$\Theta_X(\SS)$:  irreducible component of $\Psi_X(\SS)$ not of type $A_1$;  $\emptyset$ if there is no such a component\\
$\SS$ : root system \\
$\Psi_X(\LL)$:  orthogonal complement of $X$ in $\LL$\\
$\hat\LL$:  extended $\Pi$-system for an irreducible $\Pi$-system $\LL$\\

$\overline{X}$:  completion of $X$\\
$|X|$: cardinality of $X$\\
$X\simeq Y$:   isomorphic $X$ and $Y$ \\
$X\sim Y$: conjugate $X$ and $Y$\\

$(\aa,\bb)$:  entries $\frac{2(\aa,\bb)}{(\bb,\bb)}$ of Cartan matrix\\
$\<a|b\>$: $\left|(\aa|\bb)\right|$ for every $a=(\aa,-\aa), b=(\bb,-\bb)$\\
$\perp$: orthogonal (with respect to $(\cdot| \cdot)$ or to  $\<\cdot|\cdot\>$) 

\end{flushleft}


\begin{bibsection}
\begin{biblist}

\bib{Bou02}{book}{
author={Bourbaki, N.},
title={Elements of mathematics. Lie groups  Lie algebras. Chapters 4\mdash 6},
publisher={Springer},
date={2002},
}

\bib{Dyn51}{article}{
author={Dynkin, E. B.},
title={Automorphisms of simple Lie algebras},
journal={Doklady Akad Nauk SSSR N.S.},
volume={76},
date={1951},
pages={629\mdash 632},
}


\bib{Dyn52a}{article}{
author={Dynkin, E. B.},
title={Maximal subgroups of classical groups},
journal={Trudy Moskov. Mat. Obsh.},
volume={1},
date={1952},
pages={39\mdash 166},
note={ English translation in Amer. Math. Soc. Translations  (Ser. 2), Vol. 6 (1957), 111\mdash 244},
}

\bib{Dyn52}{article}{
author={Dynkin, E. B.},
title={Semisimple subalgebras of semisimple Lie algebras},
journal={Mat. Sbornik N.S.},
volume={32(72)},
date={1952},
pages={349\mdash 462},
note={ English translation in Amer. Math. Soc. Translations  (Ser. 2), Vol. 6 (1957), 111\mdash 245},
}

\bib{Hum90}{book}{
author={Humphreys, J. E.},
title={Reflection groups and Coxeter groups},
publisher={Cambridge University Press, Cambridge},
date={1990},
}

\bib{Mal44}{article}{
author={Malcev, A. I.},
title={On semisimple subgroups of Lie Groups},
journal={Izvestiya Akad. Nauk SSSR, Ser. mat.},
volume={8},
date={1944},
pages={143\mdash 174},
note={ English translation in Amer. Math. Soc. Translations, Vol. 33 (1950)},
}

\bib{McKS07}{article}{
author={McKee, J.},
author={Smyth, C.},
title={Integer symmetric matrices having all their eigenvalues in the interval~$[-2, 2]$},
journal={Journal of algebra},
volume={317(1)},
date={2007},
pages={260\mdash 290},
}

\bib{Osh07}{article}{
author={Oshima,T.},

title={A classification of subsystems of a root system},
note={arXiv: math/0611904v4 [math RT]},
date={2007},
pages={1\mdash 47}
}

\bib{Ser01}{book}{
author={Serre, J.-P.},
title={Complex semisimple Lie algebras},
publisher={Springer},
date={2001},
}

\end{biblist}
\end{bibsection}
\end{document}